\documentclass{elsart}
\usepackage{amsfonts}
\usepackage{tabmac}
\usepackage{amsthm}
\usepackage{amsmath}
\usepackage{amssymb}
\usepackage{latexsym}

\newcommand{\C}{{{\mathbb C}}}

\catcode`\@=11

\numberwithin{equation}{section}

\newtheorem*{theorem*}{Theorem}

\newtheorem{lemma}{Lemma}[section]

\newfont {\gothic}{eufm10 at 10pt}

\begin{document}

\begin{frontmatter}

\title {A Basis for the ${\bf GL_n}$ Tensor Product Algebra}
\date{\today}

\author[Yale]{Roger E. Howe}
\author[NUS]{Eng-Chye Tan}
\author[UWM]{Jeb F. Willenbring}
\address[Yale]{Department of Mathematics, Yale University, Box 208283 Yale Station, New Haven, CT 06520, USA}
\address[NUS]{Department of Mathematics, National University of Singapore, 2 Science Drive 2, Singapore 117543, Singapore}
\address[UWM]{Department of Mathematical Sciences, University of Wisconsin at Milwaukee, Milwaukee, WI 53201, USA}

\begin{abstract}
This paper focuses on the $GL_n$ tensor product algebra, which
encapsulates the decomposition of tensor products of arbitrary
irreducible representations of $GL_n$. We will describe an
explicit basis for this algebra. This construction relates
directly with the combinatorial description of
Littlewood-Richardson coefficients in terms of
Littlewood-Richardson tableaux.  Philosophically, one may view
this construction as a recasting of the Littlewood-Richardson rule
in the context of classical invariant theory.
\end{abstract}


\end{frontmatter}

\section{Introduction}

This paper continues the authors' study of the branching rules for
classical symmetric pairs by means of reciprocity algebras. In the
papers [HTW1] and [HTW2], it was shown that in some sense, the
most basic of these algebras is the tensor product algebra for the
general linear group $GL_n$. In particular, the paper [HTW1]
demonstrated that all stable branching multiplicities for
classical symmetric pairs could be written in terms of
Littlewood-Richardson coefficients, which describe  the
decomposition of  tensor products of representations of  $GL_n$.
This paper will focus on the $GL_n$ tensor product algebra. Our
main  goal is to describe an explicit basis for this algebra. As
will be seen, our basis connects directly with the combinatorial
description of Littlewood-Richardson coefficients in terms of
Littlewood-Richardson tableaux [Ful], [Pro], [Sta], [Sun].

\smallskip

For any $m \in \Bbb N$, let $B_m = A_m U_m$ be the standard Borel
subgroup of upper triangular matrices in $GL_m=GL_m(\Bbb C)$,
where $A_m$ is the diagonal torus in $GL_m$, and $U_m$ is the
maximal unipotent subgroup, i.e., the upper triangular matrices
with 1's on the diagonal. Let $\rho_m^D$ be the irreducible
representation of $GL_m$ with highest weight represented by Young
diagram $D$ as in [How].  Let $\widehat A_m^+$ be the semigroup of
dominant polynomial characters of $A_m$, i.e., those characters
which can be represented by Young diagrams.

\smallskip

Let $M_{n, (k + \ell)}$ be the space of complex $n \times (k +
\ell)$ matrices. Let $GL_n$ act by row operations (multiplication
on the left) on $M_{n, (k + \ell)}$ and $GL_k \times GL_{\ell}$ to
act on $M_{n, (k + \ell)}$ by column operations, with $GL_k$
acting on the first $k$ columns, and $GL_{\ell}$ acting on the
last $\ell$ columns.

\smallskip

Let us define a family of algebras $TA_{n, k, \ell}$, which we
shall call the {\it $GL_n$ tensor product algebra}:
$$
TA_{n, k,\ell} = P(M_{n, k + \ell})^{U_k \times U_{\ell} \times
U_n}.
$$
It is straightforward (see \S 2.1) to show that the above algebra
consists of the highest weight vectors for $GL_n$ acting on the
sum of one copy of each tensor product $\rho_n ^D \otimes \rho_n
^E$, where $D$ and $E$ are diagrams, having at most $\min \{ n, k
\}$ and $\min \{ n, \ell\}$ rows respectively. This algebra is an
$\widehat A_n^+  \times \widehat A_k^+  \times \widehat
A_{\ell}^+$-graded algebra.  In general, let $\psi^D$ denote the
$A_m$ character corresponding to the highest weight of a $GL_m$
representation $\rho_m^D$. Then, the dimension of the $\psi^F
\times \psi ^D \times \psi^E$ homogeneous subspace in $TA_{n, k,
\ell}$ equals the multiplicity of the representation $\rho_n ^F$
in the tensor product $\rho_n ^D \otimes \rho_n ^E$.  These
multiplicities are the Littlewood-Richardson coefficients.

\smallskip

Let $x_{ab}$ be the entries of a typical $n \times k$ matrix $X$,
and let $y_{ac}$ be the entries of a typical $n \times \ell$
matrix $Y$. Then
$$ Z = [X \ Y] = \left[ \begin{array}{cccccccc}
x_{11} &x_{12} &\cdots &x_{1k} &y_{11} &y_{12} &\cdots &y_{1\ell}\\
x_{21} &x_{22} &\cdots &x_{2k} &y_{21} &y_{22} &\cdots &y_{2\ell}
\\ x_{31} &x_{32} &\cdots &x_{3k} &y_{31} &y_{32} &\cdots
&y_{3\ell}\\  \vdots &\vdots &\vdots &\vdots &\vdots &\vdots
&\vdots &\vdots \\ x_{n1} &x_{n2} &\cdots &x_{nk} &y_{n1} &y_{n2}
&\cdots &y_{n\ell}
\end{array}  \right]
$$ is a typical $n \times (k + \ell)$ matrix.  Let $X_{a,b}$ denote
the upper left hand $a \times b$ submatrix of $x$'s, that
is,{\small
$$
X_{a,b} = \left [ \begin{array}{cccc} x_{11} &x_{12} &\cdots &x_{1b} \\
x_{21} &x_{22} &\cdots &x_{2b} \\ x_{31} &x_{32} &\cdots &x_{3b}
\\  \vdots &\vdots &\vdots &\vdots \\ x_{a1} &x_{a2} &\cdots
&x_{ab}
\end{array}  \right ].
$$}
Similarly, let $Y_{c,d}$ denote the upper left hand $c \times d$
submatrix of $y$'s.

\smallskip

Select partitions
$$
\begin{array}{cc}
D &= \{ d_1 \geq d_2 \geq d_3 \geq  \cdots \geq d_r
> 0 = d_{r+1}
\}, \\
E &= \{ e_1 \geq e_2 \geq e_3 \geq \cdots \geq e_s > 0 = e_{s+1}
\}, \text{and}\\
F &= \{ f_1 \geq f_2 \geq f_3 \geq \cdots \geq f_t > 0 = f_{t+1}
\},
\end{array}
$$
of depths $r \leq  k$, $s \leq \ell$ and $t \leq \min \{ n, k +
\ell \}$ respectively.  Denoting $|D| = \sum_{j \geq 0} d_j$ for
the \it size \rm of $D$, we will assume that $|D| + |E| = |F|$.  We also let $D^t$,
$E^t$ and $F^t$ denote the transpose diagrams. Thus, the
numbers $d_j$ are the lengths of the columns of $D^t$, etc.

\smallskip

Let $A=[\alpha_{ij}]$ and $B=[\beta_{i,j}]$ be $t\times r$ and $t
\times s$ matrices of indeterminates in $\alpha_{ij}$'s and
$\beta_{ij}$'s respectively.  We define $\Delta_{(D, E, F), (A,
B)}$ to be the determinant of a square $|F| \times (|D|+|E|)$
matrix as follows:  Note that the entries of this matrix are block
matrices in $x$'s or $y$'s with sizes determined by the partitions
$D$, $E$ and $F$. For instance, the block entry $\alpha_{jk}
X_{f_j,d_k}$ is an $f_j \times d_k$ matrix and likewise, the block
entry $\beta_{jk} Y_{f_j,e_k}$ is an $f_j \times e_k$ matrix.
{\small
$$
\Delta_{(D, E, F), (A, B)} =  \left|
\begin{array}{cccccccc}
\alpha_{11} X_{f_1,d_1} &\alpha_{12} X_{f_1,d_2} &\cdots
&\alpha_{1r} X_{f_1,d_r}
&\beta_{11} Y_{f_1,e_1} &\beta_{12} Y_{f_1,e_2} &\cdots &\beta_{1s} Y_{f_1,e_s}\\
\alpha_{21} X_{f_2,d_1} &\alpha_{22} X_{f_2,d_2} &\cdots
&\alpha_{2r} X_{f_2,d_r} &\beta_{21} Y_{f_2,e_1} &\beta_{22}
Y_{f_2,e_2} &\cdots
&\beta_{2s} Y_{f_2,e_s}\\
\alpha_{31} X_{f_3,d_1} &\cdots &\cdots &\vdots
&\beta_{31} Y_{f_3,e_1} &\cdots &\cdots &\vdots \\
\vdots &\vdots &\vdots &\vdots &\vdots &\vdots &\vdots &\vdots \\
\alpha_{t1} X_{f_t,d_1} &\alpha_{t2} X_{f_t,d_2} &\cdots
&\alpha_{tr} X_{f_t,d_r} &\beta_{t1} Y_{f_t,e_1} &\beta_{t2} Y_{f_t,e_2}
&\cdots &\beta_{ts} Y_{f_t,e_s}\\
\end{array} \right| \qquad
$$}
This is a polynomial function in the $x_{ab}$'s, the $y_{ac}$'s,
the $\alpha_{jk}$'s and the $\beta_{jk}$'s. Lemma \ref{lem6} shows
that for any choice of coefficient matrices $A$ and $B$, i.e.,
setting $\alpha_{jk} \in \C$ and $\beta_{jk} \in \C$, the
polynomial $\Delta_{(D, E, F), (A, B)}$ is a $GL_n \times GL_k
\times GL_{\ell}$ highest weight vector, that is, it is invariant
under $U_n \times U_k \times U_{\ell}$.

\smallskip

For generic values of $A$, we may reduce $A$ by row and column
operations to the $t \times r$ \lq \lq {\it identity}" matrix,
$J_{t,r}$, without changing the value of $\Delta_{(D, E, F), (A,
B)}$ (see \S 2.2.2). Here, $J = J_{t, r}$ is the $t \times r$
matrix with ones down the diagonal and zeroes elsewhere. Thus we
consider the polynomials $\Delta_{(D, E, F), (J,B)}$. We can
expand these as polynomials in the $\beta_{ih}$:
$$
\Delta_{(D, E, F), (J,B)} = \sum_M \Delta_{(D, E, F), M} \left
(\prod_{i, h} \beta_{ih} ^{m_{ih}} \right),
$$
where $M = [ m_{ih} ]$ is a matrix of non-negative integers,
defining the exponents to which the $\beta_{ih}$ appear in a given
monomial. Each coefficient $\Delta_{(D, E, F), M}$ is a polynomial
in $x_{ab}$ and $y_{ac}$.

\smallskip

As is well known, the Littlewood-Richardson coefficients count the
number of \it Littlewood-Richardson tableaux \rm (see \S 2.3).  We
focus on the LR\footnote{LR = Littlewood-Richardson} tableaux,
which count the multiplicities of $\rho_n ^{F^t}$ in the tensor
product $\rho_n^{D^t} \otimes \rho_n ^{E^t}$.  We will define a
mapping from Littlewood-Richardson tableaux (for appropriate fixed
$D$, $E$, $F$) to certain of the coefficients $\Delta_{(D, E, F),
M}$. More precisely, we will show how to associate to each
Littlewood-Richardson tableau $T$, a unique monomial $M(T)$ (see
\S 3.1) in the $\beta_{ih}$.

\smallskip

In essence, for fixed $D$, $E$ and $F$, a Young tableau $T$ is a
filling of a skew diagram (see \S 2.3).  The monomial $M(T)$ in
the variables $\beta_{ih}$ is determined directly from the entries
in the boxes of the skew diagram $T=F^t - D^t$.  The combinatorics
of Young tableaux leading to the definition/construction of $M(T)$
will be discussed in \S 3.1. Not surprising, one can easily and
uniquely recover $T$ from the monomial $M(T)$.

\smallskip

Now set $\Delta_{(D, E, F), M(T)}$ to be the coefficient of $M(T)$
in $\Delta_{(D, E, F), (J,B)}$. Note that $\Delta_{(D, E, F),
M(T)}$ is a polynomial in $x_{ab}$ and $y_{ac}$.  By considering a
standard coefficient in the $x$-variables of $\Delta_{(D, E, F),
M(T)}$, we have a polynomial $\delta_{T,Y}$ (see \S 3.2) in the
variables $y_{ac}$. This polynomial can also be directly related
to the entries in the boxes of $T=F^t - D^t$ (see Lemma 3.2).

\smallskip

The key result in this paper is a linear basis for the $GL_n$
tensor product algebra, whose proof is explained in \S 3:

\medskip

\begin{theorem*}
As the coefficient matrices $A$ and $B$ vary through all possible
choices of constants, the polynomials $\Delta_{(D, E, F), (A, B)}$
span the tensor product algebra $TA_{n, k, \ell}$. More precisely,
the polynomials $\Delta_{(D, E, F), M(T)}$ with the monomial
$M(T)$ associated to each Littlewood-Richardson skew tableau $T$,
form a basis for $TA_{n,k, \ell}$.
\end{theorem*}

\smallskip

The final section looks at our construction for the generators of
the $SL_4$ tensor product algebra. We correlate our monomials with
the diagrammatic description of multiplicities given by Berenstein
and Zelevinsky (see [BZ]).

\section{Preliminaries}

\subsection{$GL_n$ Tensor Product Algebra}

We begin by reviewing from [How] and [HTW2] our realization of the
tensor product algebra for $GL_n = GL_n(\C)$. All groups are
defined over $\C$ unless otherwise stated.  Let $GL_n \times GL_k$
act on the space $M_{n, k}$ of $n \times k$ complex matrices, by
the formula
\begin{equation}\label{2.1}
(g, g')(T) = (g^t)^{-1} T (g')^{-1},\qquad g \in GL_n, g' \in
GL_k, T \in M_{n,k}.
\end{equation}

Extend this action to an action by algebra automorphisms on the
ring $P(M_{n, k})$ of polynomials on $M_{n,k}$ in the standard
way. Then we have the decomposition [How]
\begin{equation}\label{2.2}
P(M_{n,k}) \simeq \sum_D \rho_n ^D \otimes \rho_k ^D,
\end{equation}
into irreducible representations for $GL_n \times GL_k$. Here
$\rho_n^D$ is the irreducible representation of $GL_n$ with
highest weight represented by Young or Ferrers diagram $D$ and the
index $D$ of summation ranges over all Young diagrams of depth not
more than $\min \{n , k \}$.

\smallskip

Let $B_n$ be the standard Borel subgroup of upper triangular
matrices in $GL_n$, and let $U_n$ be its maximal unipotent
subgroup, i.e., the upper triangular matrices with 1's on the
diagonal. Consider the algebra $P(M_{n,k})^{U_k}$ of polynomials
on $M_{n, k}$ invariant under the action of $U_k$. In terms of the
decomposition (\ref{2.2}), we may write
\begin{equation}\label{3.1}
P(M_{n,k})^{U_k} \simeq \left( \sum_D \rho_n ^D \otimes \rho_k ^D
\right)^{U_k} \simeq \sum_D \rho_n ^D \otimes (\rho_k ^D)^{U_k}.
\end{equation}

By the theory of the highest weight, the space $(\rho_k^D)^{U_k}$
is one-dimensional. Thus $P(M_{n,k})^{U_k}$ consists of one copy
of each irreducible polynomial representation of $GL_n$
corresponding to a diagram with not more than $k$ rows. Note that
the space $\rho_n ^ D \otimes (\rho_k ^D)^{U_k}$ is an eigenspace
for the diagonal torus $A_k$ in  $GL_k$, consisting of the
matrices {\small
$$
{\bf a} = \left [ \begin{array}{ccccc} a_1 &0 &0 &\cdots &0 \\
0 &a_2 &0 &\cdots &0 \\
0 &0 &a_3 &\cdots &0 \\
\vdots &\vdots &\vdots &\vdots &\vdots \\
0 &0 &0 &\cdots &a_k \end{array} \right ].
$$}
The eigenvalue of  $\bf a$ depends multiplicatively on $\bf a$:
it is a character of $A_k$. The eigencharacter of $A_k$ acting
on $\rho_n ^ D \otimes (\rho_k ^D)^{U_k}$  is equal to the highest weight of
$\rho_k^D$.  This is the character
$$
\psi^D({\bf a}) = a_1 ^{d_1} a_2 ^{d_2} \cdots a_k^{d_k},
$$
where $\{ a_j\}$ are the standard coordinates on the diagonal
torus $A_k$ of $GL_k$, and $\{ d_k \}$ are the lengths of the rows
of $D$ (see [How] for more details). Thus, the algebra
$P(M_{n,k})^{U_k}$ is a graded algebra, with grading in the
semigroup $\widehat A_k^+$ of dominant characters of $A_k$. The
graded components are irreducible representations for $GL_n$.

\smallskip

Let $M_{n, (k + \ell)}$ be the space of $n \times (k + \ell)$
matrices. We understand $GL_n$ to act by row operations
(multiplication on the left) on $M_{n, (k + \ell)}$, as in formula
(\ref{2.1}).  We understand $GL_k \times GL_{\ell}$ to act on
$M_{n, (k + \ell)}$ by column operations, with $GL_k$ acting on
the first $k$ columns, and $GL_{\ell}$ acting on the last $\ell$
columns. We may also think of $GL_k \times GL_{\ell}$ embedded in
$GL_{k+\ell}$ as block diagonal matrices:{\small
$$
(g'_1, g'_2) \sim \left [ \begin{array}{cc} g'_1 & 0 \\ 0 & g'_2
\end{array} \right ]
$$}
for $g'_1 \in GL_k$ and $g'_2 \in GL_{\ell}$. The action of $GL_k
\times GL_{\ell}$ is then just the restriction of the action
(\ref{2.1}) of $GL_{k + \ell}$.

\smallskip

With these conventions, it is more or less straightforward from
equation (\ref{3.1}) to show (see [HTW2] for details) that the
algebra
\begin{equation}\label{4.2}
TA= TA_{n, k, \ell} = P(M_{n, k + \ell})^{U_k \times U_{\ell}
\times U_n}
\end{equation}
consists of the highest weight vectors for $GL_n$ acting on the
sum of  one copy of each tensor product $\rho_n ^D \otimes \rho_n
^E$, where $D$ and $E$ are diagrams, having at most $\min \{ n, k
\}$ and $\min \{ n, \ell\}$ rows respectively. This algebra is an
$\widehat A_n^+  \times \widehat A_k^+  \times \widehat
A_{\ell}^+$-graded algebra, and the dimension of the $\psi^F
\times \psi ^D \times \psi^E$ homogeneous subspace equals the
multiplicity of the representation $\rho_n ^F$ in the tensor
product $\rho_n ^D \otimes \rho_n ^E$.  We refer to the family of
algebras $TA_{n, k, \ell}$ of equation (\ref{4.2}) as ``the $GL_n$
tensor product algebra". A more careful justification for this
terminology is given in [HTW2].

\subsection{Canonical Highest Weight Vectors}

Our goal is to construct a basis for $TA_{n, k, \ell} = TA$. To
this end, let $x_{ab}$, for $1 \leq a \leq n$, $1 \leq b \leq k$,
be standard matrix coordinates on $M_{n,k}$, and let $y_{ac}$, for
$1 \leq a \leq n$, $1 \leq c \leq \ell$, be standard matrix
coordinates on $M_{n, \ell}$. In other words, let {\small
$$
X = \left [ \begin{array}{cccc} x_{11} &x_{12} &\cdots &x_{1k} \\
x_{21} &x_{22} &\cdots &x_{2k} \\ x_{31} &x_{32} &\cdots &x_{3k} \\  \vdots &\vdots &\vdots &\vdots \\ x_{n1} &x_{n2} &\cdots &x_{nk}
\end{array}  \right ]
\quad {\rm and}\quad
Y = \left [ \begin{array}{cccc} y_{11} &y_{12} &\cdots &y_{1\ell} \\
y_{21} &y_{22} &\cdots &y_{2\ell} \\ y_{31} &y_{32} &\cdots &y_{3\ell} \\  \vdots &\vdots &\vdots &\vdots \\ y_{n1} &y_{n2} &\cdots &y_{n\ell}
\end{array}  \right ]
$$}
be typical elements of $M_{n, k}$ and $M_{n, \ell}$ respectively.

\smallskip

We combine the $x_{ab}$ and the $y_{ac}$ to get a set of coordinates on
$M_{n, (k + \ell)}$. Specifically, let
$$
Z = \left [ X \; Y \right ]
$$
be a typical element of $M_{n, (k + \ell)}$. Then $GL_k$ acts by
column operations on $X$, while $GL_{\ell}$ acts by column operations on
$Y$, and $GL_n$ acts by row operations on $X$ and $Y$ simultaneously.

Let $X_{a,b}$ denote the upper left hand $a \times b$ submatrix of
$X$. That is,{\small
\begin{equation}\label{5.2}
X_{a,b} = \left [ \begin{array}{cccc} x_{11} &x_{12} &\cdots &x_{1b} \\
x_{21} &x_{22} &\cdots &x_{2b} \\ x_{31} &x_{32} &\cdots &x_{3b} \\  \vdots &\vdots &\vdots &\vdots \\ x_{a1} &x_{a2} &\cdots &x_{ab}
\end{array}  \right ].
\end{equation} }
Similarly, let $Y_{c,d}$ denote the upper left hand $c \times d$
submatrix of $Y$.

\smallskip

Following standard yoga, we think of a diagram alternatively as a partition,
and/or as a decreasing sequence of non-negative integers, with two sequences
whose positive  elements coincide being considered equivalent. Select
partitions

$$
D = \{ d_1 \geq d_2 \geq d_3 \geq  \cdots \geq d_r > 0 = d_{r+1} \},
$$
\begin{equation}\label{5.3}
E = \{ e_1 \geq e_2 \geq e_3 \geq \cdots \geq e_s > 0 = e_{s+1}
\},
\end{equation}
and
$$
F = \{ f_1 \geq f_2 \geq f_3 \geq \cdots \geq f_t > 0 = f_{t+1} \},
$$
of depths $r \leq  k$, $s \leq \ell$ and $t \leq \min \{ n, k + \ell \}$ respectively.  Let $|D| = \sum_{j \geq 0} d_j$ be the \it size \rm of
$D$. We will assume that $|D| + |E| = |F|$.

Let{\small
\begin{equation}\label{5.4}
A = \left [ \begin{array}{cccc} \alpha_{11} &\alpha_{12} &\cdots &\alpha_{1r} \\
\alpha_{21} &\alpha_{22} &\cdots &\alpha_{2r} \\ \alpha_{31}
&\alpha_{32} &\cdots &\alpha_{3r} \\  \vdots &\vdots &\vdots
&\vdots \\ \alpha_{t1} &\alpha_{t2} &\cdots &\alpha_{tr}
\end{array}  \right ].
\end{equation} }
be a $t \times r$ matrix of indeterminates $\alpha_{jk}$.
Similarly, let {\small
$$
B = \left [ \begin{array}{cccc} \beta_{11} &\beta_{12} &\cdots &\beta_{1s} \\
\beta_{21} &\beta_{22} &\cdots &\beta_{2s} \\ \beta_{31} &\beta_{32} &\cdots &\beta_{3s} \\
\vdots &\vdots &\vdots &\vdots \\ \beta_{t1} &\beta_{t2} &\cdots &\beta_{ts}
\end{array}  \right ].$$ }
be a $t \times s$ matrix of other indeterminates $\beta_{jk}$.

\smallskip

Define a matrix $\widetilde X_{(D, F, A)} = \widetilde X$ as
follows. $\widetilde X$ is a $|F| \times |D|$ matrix, which is
partitioned into submatrices according to the partition $F$ of the
rows of $\widetilde X$, and the partition $D$ of the columns of
$\widetilde X$. Thus $\widetilde X$ consists of submatrices
$\widetilde X_{j,k}$ for $1 \leq j \leq t$ and $1 \leq k \leq r$:
{\small
\begin{equation}\label{5.6}
\widetilde X \! =\!  \left [ \begin{array}{cccc}
\widetilde X_{1,1} &\widetilde X_{1,2} &\cdots &\widetilde X_{1,r} \\
\widetilde X_{2,1} &\widetilde X_{2,2} &\cdots &\widetilde X_{2,r} \\
\widetilde X_{3,1} &\cdots &\cdots &\vdots \\
\vdots &\vdots &\vdots &\vdots \\
\widetilde X_{t,1} &\widetilde X_{t,2} &\cdots &\widetilde X_{t,r} \end{array} \right ]
\! =\! \left [ \begin{array}{cccc}
\alpha_{11} X_{f_1,d_1} &\alpha_{12} X_{f_1,d_2} &\cdots &\alpha_{1r} X_{f_1,d_r} \\
\alpha_{21} X_{f_2,d_1} &\alpha_{22} X_{f_2,d_2} &\cdots &\alpha_{2r} X_{f_2,d_r} \\
\alpha_{31} X_{f_3,d_1} &\cdots &\cdots &\vdots \\
\vdots &\vdots &\vdots &\vdots \\
\alpha_{t1} X_{f_t,d_1} &\alpha_{t2} X_{f_t,d_2} &\cdots &\alpha_{tr}
X_{f_t,d_r} \end{array} \right ]\qquad
\end{equation}  }
The submatrix $\widetilde X_{j,k}$ is an $f_j \times d_k$ matrix. Specifically,
\begin{equation}\label{5.7}
\widetilde X_{j,k} = \alpha_{jk} X_{f_j,d_k}
\end{equation}
with $\alpha_{jk}$ as in formula (\ref{5.4}), and $X_{f_j, d_k}$ as
in formula (\ref{5.2}).

\smallskip

We also define a matrix $\widetilde{Y}_{(E, F, B)} = \widetilde{Y_{\ }}$ in parallel
fashion, but using the matrix $Y$ instead of $X$, the partition $E$ instead
of $D$, and the matrix $B$ instead of $A$. Finally, we define a matrix
\begin{equation}\label{5.8}
\widetilde{Z_{\  }}= \left [  \widetilde X \ \widetilde{Y_{\ }} \right ].
\end{equation}
We note that $\widetilde{Z_{\ }}$ is an
$|F| \times (|D| + |E|) = |F| \times |F|$ matrix. It is a square matrix, so
we may take its determinant:
\begin{equation}\label{5.9}
\Delta_{(D, E, F), (A, B)} = \Delta_{(D, E, F), (A, B)} (X, Y) =
\det \widetilde{Z_{\ }}.
\end{equation}
This is a polynomial function in the $x_{ab}$'s, the $y_{ac}$'s, the $\alpha_{jk}$'s and  the $\beta_{jk}$'s.

\medskip

\begin{lemma}\label{lem6}
\begin{enumerate}
\item[(a)] For any choice of coefficient matrices $A$ and $B$, the
polynomial \linebreak $\Delta_{(D, E, F), (A, B)}$ is a $GL_n
\times GL_k \times GL_{\ell}$ highest weight vector, that is, it
is invariant under $U_n \times U_k \times U_{\ell}$. \item[(b)]
$\Delta_{(D, E, F), (A, B)}$ is a weight vector for the product
$A_n \times A_k \times A_{\ell}$ of diagonal tori of the groups
$GL_n$, $GL_k$ and $GL_{\ell}$. The eigencharacter of $A_n \times
A_k \times A_{\ell}$ defined by $\Delta_{(D, E, F), (A, B)}$ is
$\psi^{F^t} \times \psi^{D^t} \times \psi^{E^t}$.
\end{enumerate}
\end{lemma}

\smallskip

\noindent {\bf Proof:} If we can show that $\Delta=\Delta_{(D, E,
F), (A, B)}$ is annihilated by the infinitesimal generators of the
groups $U_n$, $U_k$ and $U_{\ell}$, the lemma will follow. The
infinitesimal generators of $U_n$ act on $\Delta$ by row
operations. Specifically, we have a basis of the infinitesimal
generators for $U_n$, consisting of operators
$$
E_{ad} = \sum_{b = 1}^k x_{ab} {\partial \over \partial x_{db}} +
\sum_{c = 1}^{\ell} { y_{ac} {\partial \over \partial y_{dc}}},
$$
for $1 \leq a < d \leq n$. Thus, $E_{ad} \Delta$ is a sum of
terms, each of which is a determinant of a matrix in which one of
the rows of $\Delta$  consisting of entries $\alpha_{ij} x_{db}$
or $\beta_{i h} y_{dc}$ gets replaced by a row with entries
$\alpha_{ij} x_{ab}$ or $\beta_{ih} y_{ac}$. Since the block
matrices $X_{f_j, d_k}$ of formula (\ref{5.2}) contain all rows of
$X$ up to the $f_j$-th, we see that each term in the sum defining
$E_{ad} \Delta$ will be the determinant of a matrix with two
repeated rows, and therefore vanishes identically. We conclude
that $\Delta$ is invariant under $U_n$. A similar argument, but
using column operations, shows that $\Delta$ is invariant under
$U_k$ and $U_\ell$. This proves part (a) of the Lemma.  Part (b)
of the lemma is a straightforward calculation. $\qquad \square$

\subsubsection{Row and Column Operations on $\Delta_{(D, E, F), (A, B)}$}

The polynomials $\Delta_{(D, E, F), (A, B)}$ give us a large
collection of elements of $TA_{n, k, \ell}$. We want to show that
they in fact span $TA_{n, k, \ell}$. More precisely, we want to
show that, by making appropriate choices of the coefficient
matrices $A$ and $B$, we can extract from the polynomials
$\Delta_{(D, E, F), (A, B)}$ a basis for the $\psi^{F^t} \times
\psi^{D^t} \times \psi^{E^t}$-eigenspace of $TA_{n, k, \ell}$.

\smallskip

Although the $\Delta_{(D, E, F), (A, B)}$ were constructed as
polynomials in the $x_{ab}$'s and $y_{ac}$'s, they in fact also
depend polynomially on the auxiliary variables $\alpha_{ij}$ and
$\beta_{ih}$. Furthermore, these auxiliary variables are organized
into matrices, and thus support an action of a product of general
linear groups. Precisely, there is a natural action of $GL_t
\times GL_r$ on the $\alpha_{ij}$ and an action of $GL_t \times
GL_s$ on the $\beta_{ih}$. (Here, $r$, $s$ and $t$ are the depths
of the partitions $D$, $E$ and $F$ respectively (see equation
(\ref{5.3})). The $\Delta_{(D, E, F), (A, B)}$ have a very nice
property with respect to these actions.

\smallskip

Let $\widetilde U_r$ denote the group of strict left-to-right
column operations on the matrix $A$. (The tilde on $\widetilde
U_r$ is simply to distinguish the action on the auxiliary
variables $\alpha_{ij}$ from the action on the variables $x_{ij}$
and $y_{ij}$.)  That is, $\widetilde U_r$ allows one to add a
multiple of a column of the matrix $A$ to any column strictly to
its right, and to combine such operations. It is a maximal
unipotent subgroup of $GL_r$. Similarly, let $\widetilde U_s$ be
the strict left-to-right column operations on $B$. Finally, let
$\widetilde U_t$ be the diagonal group of strict top-to-bottom row
operations on $A$ and $B$ simultaneously. This is generated by the
operations which add a multiple of one row of $A$ to a lower row
of $A$, and at the same time, do the same thing to $B$. Let
$\widetilde A_r$ be the diagonal torus of $GL_r$, which multiplies
each column of $A$ by some scalar. Define $\widetilde A_s$ and
$\widetilde A_t$ analogously.

\medskip

\begin{lemma}\label{lem7}
As functions of the $\alpha_{ij}$ and the $\beta_{ih}$, the
polynomials $\Delta_{(D, E, F), (A, B)}$ are invariant under the
action of $\widetilde U_t \times \widetilde U_r \times \widetilde
U_s$. They are eigenvectors for the product $\widetilde A_t \times
\widetilde A_r \times \widetilde A_s$, with eigencharacter
$\psi^{F} \times \psi^{D} \times \psi^{E}$.
\end{lemma}

\smallskip

\noindent {\bf Proof:} Consider the column operations which add
$\eta$ times the $i$-th column of $A$ to the $h$-th column, where
$i < h$. That is, this operation replaces the entries
$\alpha_{jh}$, $1 \leq j \leq t$, with entries $\alpha_{jh} + \eta
\alpha_{ji}$. Here $\eta$ is a complex number. From the form of
the entries of $\widetilde X$ (see formula (\ref{5.7})), and the
fact that the numbers $d_i$ decrease as $i$ increases, we can see
that the effect of this operation on $\widetilde X$ can also be
achieved by means of column operations directly on $\widetilde X$.
Precisely, for each of the first $d_h$ columns of the $i$-th block
of columns of $\widetilde X$, we should add $\eta$ times a column
to the corresponding column of the $h$-th block of columns. Since
these operations are all strict left-to-right column operation is,
they do not change the value of $\Delta_{(D, E, F), (A, B)}$. The
invariance statement of the lemma follows. The statement about the
behavior under diagonal operations is again a straightforward
computation.  $\qquad \square$

\subsubsection{The Reduced Polynomials $\Delta_{(D, E, F), (J,B)}$}

By the Lemma \ref{lem7}, we see that, for generic values of $A$,
we may reduce $A$ to the $t \times r$ \lq \lq {\it identity}"
matrix, $J_{t,r}$, without changing the value of $\Delta_{(D, E,
F), (A, B)}$. Here, $J = J_{t, r}$ is the $t \times r$ matrix with
ones down the diagonal and zeroes elsewhere. In other words, if $A
=LV$, where $L$ is a lower triangular $t \times r$ matrix, and $V$
is an upper triangular $r \times r$ unipotent matrix, then we have
\begin{equation}\label{8.1}
\Delta_{(D, E, F), (A, B)} = \left( \prod_{i = 1}^t
\alpha_{ii}^{f_i} \right) \Delta_{(D, E, F), (J,B)}.
\end{equation}
In particular, the polynomials of the form $\Delta_{(D, E, F),
(J,B)}$ will span the same space as the collection of all
$\Delta_{(D, E, F), (A, B)}$.

\smallskip

Thus we consider  the polynomials $\Delta_{(D, E, F), (J,B)}$. We
can expand these as polynomials in the $\beta_{ih}$:
\begin{equation}\label{8.2}
\Delta_{(D, E, F), (J,B)} = \sum_M \Delta_{(D, E, F), M} \left
(\prod_{i, h} \beta_{ih} ^{m_{ih}} \right),
\end{equation}
where $M = [ m_{ih} ]$ is a matrix of non-negative integers,
defining the exponents to which the $\beta_{ih}$ appear in a given
monomial. Each coefficient $\Delta_{(D, E, F), M}$ is a polynomial
in $x_{ab}$ and $y_{ac}$.

\subsection{The Structure of Littlewood-Richardson Tableaux}

As is well known, the Littlewood-Richardson (which we abbreviate
using LR) coefficients count the number of \it
Littlewood-Richardson tableaux \rm [CGR], [RW], [Ful], [Sta]. With
fixed and appropriate $D$, $E$ and $F$, we focus on LR tableaux
which capture the multiplicity of $\rho_n^{F^t}$ in the tensor
product $\rho_n^{D^t} \otimes \rho_n^{E^t}$ (see [Ful]). We will
define a mapping from Littlewood-Richardson tableaux (for
appropriate fixed $D$, $E$ and $F$) to certain of the coefficients
$\Delta_{(D, E, F), M}$. More precisely, we will show how to
associate to each Littlewood-Richardson tableau $T$, a monomial
$M(T)$ in the $\beta_{ih}$.

\smallskip

A Young tableau $T$ is a filling of a skew diagram. Fix $D$, $E$
and $F$. This monomial $M(T)$ in the variables $\beta_{ij}$ is
directly related to the entries of the skew diagram $T=F^t - D^t$.
The related combinatorics of Young tableaux leading to the
definition/construction of $M(T)$ will be discussed in \S 3.1. Not
surprising, one could easily and uniquely recover $T$ from the
monomial $M(T)$.  We will also associate to $T$ the polynomial
$\Delta_{(D, E, F), M(T)}$, i.e., the coefficient of $M(T)$ in the
expansion (\ref{8.2}) of $\Delta_{(D, E, F),(J,B)}$.

\smallskip

Note that $\Delta_{(D, E, F), M(T)}$ is a polynomial in $x_{ab}$
and $y_{ac}$.  By considering a standard coefficient in the
$x$-variables of $\Delta_{(D, E, F), M(T)}$, we have a polynomial
$\delta_{T,Y}$ (see \S 3.3) in the variables $y_{ac}$. This
polynomial can also be directly related to the fillings of $T=F^t
- D^t$ (see Lemma 3.1 and Lemma 3.2).  Lemma 3.2 is, in fact, the
key to the proof of our linear basis result.

\smallskip

Constructing the correspondence as proposed above will involve
some study of the structure of LR tableaux. We recall that a
Littlewood-Richardson tableau (= LR tableau) is a skew diagram
filled with positive integers 1 through $k$. The criteria for such
a tableau $T$ to be LR may be stated as follows [CGR], [Ful]:
\begin{enumerate}

\item[LR1:] $T$ is semistandard, which means that the numbers in
each row of $T$ weakly increase from left to right, and the
numbers in each column strictly increase from top to bottom.

\item[LR2:] For every pair of positive numbers $m$ and $p$, the
number of times the number $m$ occurs in the first $p$ rows of $T$
is not larger than the number of times that $m-1$ appears in the
first $p-1$ rows of the tableau. This condition is interpreted as
being vacuous when $m =1$.
\end{enumerate}

\smallskip

\noindent {\bf Example:}  As far as possible, we will try to use the same example throughout this paper, and this example has to be complicated enough for us to point out several features which are important. Our example will be as follows:
$$
D=(3,3,2,1,1),i.e., r=5 \qquad {\rm and} \qquad
D^t =\tableau[s]{ & & & & \\ & & \\ &  }
$$

$$
E=(3,3,2,1),i.e., s=4 \qquad {\rm and} \qquad E^t =\tableau[s]{   & & & \\ & & \\
&  }
$$

$$
F=(5,5,4,3,1,1),i.e., t=6 \qquad {\rm and} \qquad
$$

$$
T =\tableau[s]{ & & & & &\tf{1} \\ & & &\tf{1} \\ & &\tf{1} &\tf{2} \\ \tf{1} &\tf{2} &\tf{3} \\ \tf{2} &\tf{3} },\qquad
T_1 =\tableau[s]{ & & & & &\tf{1} \\ & & &\tf{1} \\ & &\tf{1} &\tf{2} \\ \tf{1} &\tf{2} &\tf{2} \\ \tf{3} &\tf{3} }
$$

\smallskip

$$
T_2 =\tableau[s]{ & & & & &\tf{1} \\ & & &\tf{1} \\ & &\tf{2} &\tf{2} \\ \tf{1} &\tf{1} &\tf{3} \\ \tf{2} &\tf{3} }    \qquad {\rm and} \qquad
T_3 =\tableau[s]{ & & & & &\tf{1} \\ & & &\tf{2} \\ & &\tf{1} &\tf{3} \\ \tf{1} &\tf{1} &\tf{2} \\ \tf{2} &\tf{3} }
$$
are four tableaux for $F^t - D^t$.  Observe that the skew tableaux
$T$, $T_1$, $T_2$ and $T_3$ are fillings of the tableau $F^t -
D^t$, i.e., the difference between diagrams $F^t$ and $D^t$, as
depicted by the cells whose boundaries are highlighted.  (For
those familiar with the combinatorics of Young diagrams, note that
each of the tableau must have weight $E^t$.)  We can verify that
the fillings are the only possible fillings satisfying conditions
(LR1) and (LR2), and thus the LR coefficient $c_{D^t,E^t}^{F^t}
=4$. This of course means that the irreducible representation with
highest weight parameterized by $F^t$ appears four times in the
tensor product of irreducible representations parameterized by
$D^t$ and $E^t$.

\subsubsection{Banal Tableau}

By a \it banal \rm LR tableau $BY$, we will mean an ordinary
(i.e., not skew) diagram $Y$ in which the cells of the $i$-th row
are all filled with the number $i$. Equivalently, a column of
length $k$ is filled with the numbers from 1 to $k$, in increasing
order as you move down the column. We will describe a content
preserving mapping between the cells of a general LR tableau and a
banal LR tableau. By \lq \lq {\it content preserving}", we mean
that each cell of one tableau is mapped to a cell of the other
tableau with the same value (i.e., same number inside).

\smallskip

\noindent {\bf Example:} In the above example where
$$
E=(3,3,2,1)\qquad {\rm we\ have\ banal\ tableau\ for\ } E^t =BE^t=
\tableau[s]{  1& 1& 1& 1 \\ 2 &2 &2 \\ 3 &3 }
$$

\subsubsection{Peeling a Tableau $T$}

In any diagram or skew-diagram, a cell which lies above and to the
right of a given cell will be said to lie \it northeast \rm of the
given cell.  Note that northeast also includes things in the same
row (resp. column), but to the right (resp. above).

\smallskip

Consider a general LR tableau $T$. Let $\ell_0$ be the largest number which
occurs in any cell of $T$. For each number $h$, $1 \leq h \leq \ell_0$, let
$C_1(h)$ be the cell of $T$ which contains $h$, and which lies farthest to
the northeast among all cells containing $h$.

\smallskip

\noindent{\bf Example:} For our example in \S 2.3, $\ell_0=3$, and
we have highlighted $C_1(1)$, $C_1(2)$ and $C_1(3)$ for both
tableaux $T$ and $T_3$: Note that $F=(5,5,4,3,1,1)$ and
$$
T =\tableau[s]{ & & & & &\tf{1} \\ & & &1 \\ & &1 &\tf{2} \\ 1 &2 &\tf{3} \\ 2 &3 }    \qquad {\rm and} \qquad
T_3 =\tableau[s]{ & & & & &\tf{1} \\ & & &\tf{2} \\ & &1 &\tf{3} \\ 1 &1 &2 \\ 2 &3 }
$$

\smallskip

We should check that $C_1 (h)$ is well-defined. Start with $C_1(\ell_0)$.
Let $C$ be a cell containing $\ell_0$, such that there is no cell
containing $\ell_0$ and lying to the northeast of $C$. Then $C$ must lie on
the end of a row, since rows are weakly increasing to the right, and
$\ell_0$ is the largest possible entry a cell can have.

\smallskip

Also, no cell in a row above the row of $C$ can contain $\ell_0$. For, if
$C'$ is such a cell, then also the cell $C''$ at the right end of the row of
$C'$, must contain $\ell_0$, since the entries must weakly increase to the
right, hence be at least equal to $\ell_0$, but also cannot be larger. But
the end of each row of a skew diagram is to the northeast of the ends of lower
rows.  Thus, by choice of $C$, the cell $C''$ cannot exist, and therefore
also, $C'$ cannot exist. Therefore, we must have that $C = C_1(\ell_0)$, so
$C_1(\ell_0)$ is indeed well-defined.

\smallskip

Let $T^+$ be the tableau consisting of the rows of $T$ above the
row of $C_1(\ell_0)$. By the argument of the previous paragraph,
the entries of $T^+$ are all less than $\ell_0$. On the other
hand, the condition (LR2) for LR tableaux guarantees that some
cells of $T^+$ do contain $\ell_0 - 1$, which is thus the largest
entry in $T^+$. Therefore, $C_1(\ell_0 - 1)$ will be contained in
$T^+$, and the argument just given implies that $C_1(\ell_0 -1)$
does indeed exist. Continuing by downward induction, we see that
$C_1(h)$ is well-defined for all $h$. Furthermore, as byproducts
of the existence argument, we have found that each $C_1(h)$ lies
on the right end of its row, and that $C_1(h-1)$ lies to the
northeast of $C_1(h)$, more precisely, $C_1(h-1)$ lies strictly
above and weakly to the right of $C_1 (h)$.  In our example $T_3$,
the boxes $C_1(2)$ and $C_1(3)$ lie in the same column.  Thus, the
cells $C_1(h)$ constitute what is sometimes called a \it vertical
skew strip. \rm  In our example above, this vertical strip has
been highlighted.

\smallskip

We now remove from $T$ the cells $C_1(h)$ for $1 \leq h \leq
\ell_0$, and reassemble them into a column of length $\ell_0$, in
consecutive order from top to bottom. This column will be the
first column of the banal tableau to which we will map $T$. We
will denote this banal tableau by $BT$.

\smallskip

\noindent{\bf Example:} For our example in \S 2.3, $\ell_0=3$, and
we have highlighted $C_1(1)$, $C_1(2)$ and $C_1(3)$: Here
$F=(5,5,4,3,1,1)$ and
$$
T =\tableau[s]{ & & & & &\tf{1} \\ & & &1 \\ & &1 &\tf{2} \\ 1 &2 &\tf{3} \\ 2 &3 }  \qquad \longrightarrow \qquad
\tableau[s]{ 1 \\ 2 \\ 3} \quad {\rm to\ be\ attached\ to\ } BT'.
$$

\medskip

\begin{lemma}\label{lem13}
Let $T'$ denote the configuration which is left after removal of
the cells $C_1(h)$ from $T$. We claim that $T'$ is an LR tableau.
\end{lemma}

\smallskip

\noindent{\bf Example:} For our example in \S 2.3, $\ell_0=3$,
after removal of the vertical strip formed by $C_1(1)$, $C_1(2)$
and $C_1(3)$: Recall that $F=(5,5,4,3,1,1)$ and
$$
T' =\tableau[s]{ & & & & \\ & & &1 \\ & &1 \\ 1 &2 \\ 2 &3 }
\qquad {\rm and} \qquad
T_3' =\tableau[s]{ & & & & \\ & & \\ & &1 \\ 1 &1 &2 \\ 2 &3  }
$$
are again LR tableaux.

\smallskip

\noindent {\bf Proof:} First, we argue that the cells of $T'$ form
a skew diagram. We have seen that each cell $C_1(h)$ lies at the
end of its row. Furthermore, $C_1(\ell_0)$ must lie also at the
bottom of its column, since it contains the largest number in $T$,
and columns are strictly increasing as you go down.  Hence  $T -
C_1({\ell_0})$ (i.e., removal of the cell $C_1(\ell_0)$ from $T$)
is a skew diagram.

\smallskip

Suppose that $T$ with the cells $C_1(h)$ removed for $h > g$ is a skew
diagram. Call it $T_{(g)}$ (not to be confused with the tableaux $T_1$, $T_2$ and $T_3$ in our example).  Now let us remove the cell $C_1 (g)$. As with all
the cells $C_1(h)$, the cell  $C_1 (g)$ lies at the end of its row, even in $T$, hence
\underbar{a fortiori} in $T_{(g)}$. If there is a cell in $T_{(g)}$ which lies below $C_1(g)$, the number in this cell must be greater than $g$. Say it is $g' > g$. Since $C_1(g)$ is the farthest northeast occurrence of $g$, the cell immediately below mus

t be northeast of $C_1(g')$, contradicting the definition of $C_1 (g')$. Hence, there is no cell in $T_{(g)}$ lying directly below $C_1 (g)$, so $T_{(g)} - C_1 (g)$ is a skew diagram. By downward induction from $\ell_0$, we conclude that the cells of $T'$

 do form a diagram, as desired.

\smallskip

Since $T'$ is a subdiagram of $T$, it is clear by reference to
condition (LR1) that $T'$ is semistandard.  To show that $T'$ is
an $LR$ tableau, we have still to verify that condition (LR2) is
valid. Hence consider positive integers $m$ and $p$. We must argue
that the number of times that $m$ occurs in the first $p$ rows of
$T'$ is at least as large as the number of times that $m+1$ occurs
in the first $p+1$ rows of $T'$. If $m+1$ does not occur in the
first $p+1$ rows of $T'$, then this is clearly true. If however
$m+1$ does occur in the first $p+1$ rows of $T'$, then it also
certainly occurs in the first $p+1$ rows of $T$, and therefore the
cell $C_1(m+1)$ must have belonged to the first $p+1$ rows of $T$.
It follows that the number of occurrences of $m+1$ in the first
$p+1$ rows of $T'$ is one less than the number of occurrences of
$m+1$ in the first $p+1$ rows of $T$. On the other hand, the
number of occurrences of $m$ in the first $p$ rows of $T'$ is at
least one less than  the number of occurrences of $m$ in the first
$p$ rows of $T$, and this in turn is at least equal to the number
of occurrences of $m+1$ in the first $p+1$ rows of $T$, since $T$
is LR. It follows that condition (LR2) is also satisfied by $T'$.
Therefore, $T'$ is LR. $\qquad \square$

\subsubsection{Standard Peeling of a Tableau $T$}

We can now finish our description of the mapping from $T$ to the
banal tableau $BT$. We have removed from $T$ the boxes $C_1(h)$,
for $1 \leq h \leq \ell_0$, and reassembled them into the first
column of the prospective $BT$. This has left us with a
configuration $T'$ of cells, and we have shown that $T'$ is an LR
tableau. By induction on the number of boxes in $T$, we may assume
that the process we have outlined above allows us to define the
banal tableau $BT'$, and the mapping from $T'$ to $BT'$. To get
$BT$, we append to the left of $BT'$, the column of length
$\ell_0$ which we have just constructed.  More precisely, the
first column of $BT$ is this column, and for $i > 1$, the $i$-th
column of $BT$ is the $(i-1)$-th column of $BT'$. Since the
largest entry in any cell of $T'$ is at most $\ell_0$, the columns
of $BT'$ cannot be longer than $\ell_0$, and therefore the shape
of $BT'$ is a diagram, and the entries in its $i$-th row will
clearly all be equal to $i$. Thus, $BT$ is a banal tableau. As for
the mapping from $T$ to $BT$, we take the union of the mappings
from $T'$ to $BT'$ and the mapping of the cells $C_1(h)$ to the
first column of $BT$. This evidently is content preserving. Hence,
we have constructed the desired mapping.

\smallskip

Since the construction of $BT$ from $T$ involves successive removal of
vertical skew strips from $T$, as described in \S 2.3.2, we will refer to the
mapping from $T$ to $BT$ constructed above, in
\S 2.3.2 to \S 2.3.3, as the {\it standard peeling} of $T$.

\smallskip

\noindent{\bf Example:} The following gives an illustration of the
standard peeling of $T$ for our example in \S 2.3. The vertical
strips are being highlighted at each stage of peeling: {\small
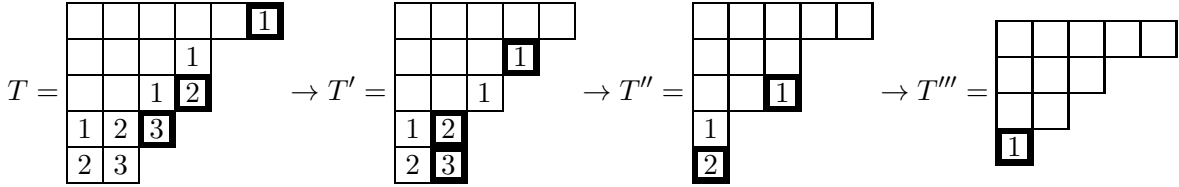
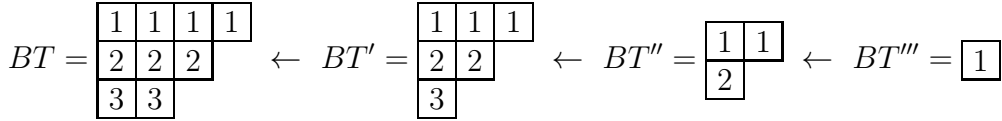
$$
T =\tableau[s]{ & & & & &\tf{1} \\ & & &1 \\ & &1 &\tf{2} \\ 1 &2 &\tf{3} \\ 2 &3 }   \rightarrow
T' =\tableau[s]{ & & & & \\ & & &\tf{1} \\ & &1 \\ 1 &\tf{2} \\ 2 &\tf{3} }  \rightarrow
T'' =\tableau[s]{ & & & & \\ & & \\ & &\tf{1} \\ 1 \\ \tf{2} }
 \rightarrow
T''' =\tableau[s]{ & & & & \\ & & \\ & \\ \tf{1} }
$$
}

$$
BT=\tableau[s]{ 1 &1 &1 &1 \\ 2 &2 &2 \\ 3 &3 }
\ \leftarrow \
BT'=\tableau[s]{ 1 &1 &1 \\ 2 &2 \\ 3 }
\ \leftarrow \
BT''=\tableau[s]{ 1 &1 \\ 2 }
\ \leftarrow \
BT'''=\tableau[s]{ 1 }
$$

\section{Proof of Theorem}

\subsection{The Monomial $M(T)$, in the
Variables $\beta_{ih}$, Associated to $T$}

With the discussions of \S 2.2 and \S 2.3 in mind, we can now describe
how to associate to an LR tableau $T$, a monomial $M(T)$ in the
auxiliary variables $\beta_{ih}$.

\smallskip

The LR tableau lives in a skew-diagram, but we assume that this
skew-diagram is embedded in a larger diagram, associated to the
partition $F^t$.  The skew diagram will be the difference $F^t -
D^t$, between $F^t$ and a smaller diagram $D^t$. The LR tableau
$T$ will result from a filling of $F^t - D^t$ by the entries from
a banal tableau $BE^t$. The inverse of this filling will be the
standard peeling from $F^t - D^t$ to  $E^t = B(F^t- D^t)$. The
partitions $D$, $E$ and $F$ are as in (\ref{5.3}).

\smallskip

Consider the standard peeling of the LR tableau $T$, or rather,
its inverse, the filling of the shape of $F^t - D^t$ by the
contents of $BE^t$. The elements of the $h$-th column of $BE^t$
will be distributed among various columns of $F^t - D^t$. We
number the columns of each diagram consecutively, from left to
right. Thus, the $h$-th column of $E^t$ has length $e_h$, and
there are $s$ columns all together (so $s$ is the length of the
first row of $E^t$). Suppose that $m_{ih}$ of the elements from
the $h$-th column of $BE^t$ get put into the $i$-th column of $F^t
- D^t$. Then we associate to the tableau $T$ the monomial
\begin{equation}\label{15.1}
M(T) = \prod_{i,h} \beta_{ih}^{m_{ih}}.
\end{equation}

As indicated above in \S 2.2.2, we now associate to $T$ the
polynomial $\Delta_{(D, E, F), M(T)}$  in the  $x_{ab}$ and
$y_{ac}$, which is the coefficient of $M(T)$ in the expansion
(\ref{8.2}):
$$
\Delta_{(D, E, F), (J,B)} = \sum_M \Delta_{(D, E, F), M} \left
(\prod_{i, h} \beta_{ih} ^{m_{ih}} \right).
$$

\smallskip

\noindent {\bf Example:} In our example from \S 2.3, we have
$$
\Delta_{(D,E,F), (J, B)} = {\small \left| \begin{array}{ccccccccc}
X_{5,3} &0 &0 &0 &0
&\beta_{11}Y_{5,3} &\beta_{12}Y_{5,3} &\beta_{13}Y_{5,2} &\beta_{14}Y_{5,1} \\
0      &X_{5,3}&0 &0 &0
&\beta_{21}Y_{5,3} &\beta_{22}Y_{5,3} &\beta_{23}Y_{5,2} &\beta_{24}Y_{5,1} \\
0 &0 &X_{4,2} &0 &0
&\beta_{31}Y_{4,3} &\beta_{32}Y_{4,3} &\beta_{33}Y_{4,2} &\beta_{34}Y_{4,1} \\
0 &0 &0 &X_{3,1} &0
&\beta_{41}Y_{3,3} &\beta_{42}Y_{3,3} &\beta_{43}Y_{3,2} &\beta_{44}Y_{3,1} \\
0 &0 &0 &0 &X_{1,1}
&\beta_{51}Y_{1,3} &\beta_{52}Y_{1,3} &\beta_{53}Y_{1,2} &\beta_{54}Y_{1,1} \\
0 &0 &0 &0 &0
&\beta_{61}Y_{1,3} &\beta_{62}Y_{1,3} &\beta_{63}Y_{1,2} &\beta_{64}Y_{1,1} \end{array} \right|}
$$
Note that $F= (5,5,4,3,1,1)$ gives the heights of the blocks, from top to bottom. The widths of the blocks (from left to right) are $D=(3,3,2,1,1)$ for the $x$'s and $E=(3,3,2,1)$ for the $y$'s.

\smallskip

\noindent {\bf Example:} For the tableau $T$, the following matrix
gives the entries $m_{ih}$ of the elements from the $h$-th column
of $BE^t$ get put into the $i$-th column of $F^t - D^t$: {\small
$$
\left[ \begin{array}{cccc} 0 &0 &1 &1 \\ 0 &2 &0 &0 \\ 1 &0 &1 &0 \\
1 &1 &0 &0 \\ 0 &0 &0 &0 \\ 1 &0 &0 &0 \end{array}\right]
$$}
and so $M(T)= \beta_{31}\beta_{41}\beta_{61}\beta_{22}^2\beta_{42}\beta_{13}\beta_{33}
\beta_{14}.$

\smallskip
Compare with that of tableau $T_1$: {\small
$$
\left[ \begin{array}{cccc} 0 &1 &0 &1 \\ 1 &0 &1 &0 \\ 0 &1 &1 &0 \\
1 &1 &0 &0 \\ 0 &0 &0 &0 \\ 1 &0 &0 &0 \end{array}\right]
$$ }
with
$M(T_1)= \beta_{21}\beta_{41}\beta_{61}\beta_{12}\beta_{32}\beta_{42}
\beta_{23}\beta_{33}\beta_{14}.$

\smallskip

\noindent There is a combinatorial formula for $\Delta_{(D,E,F),
M(T)}$ (see Lemma 4.2.2 in [HL]). It is a partial expansion of a
certain determinant. For our example, there is an even simpler
expression in terms of sums of determinants. The following two
examples will illustrate this.

Principally, given a skew tableau $T$ with $M(T)=\Pi_{i,h}
\beta_{ih}^{m_{ih}}$ defined by the filling of $T$ as in (3.1), we
can look at coefficients $\Pi_{i,h} \beta_{ih}^{n_{ih}} $ of
$\Delta_{(D,E,F),(J,B)}$ where the $t \times s$ matrices of
exponents $N=[n_{ih}]$ satisfy the following properties:
\begin{enumerate}
\item[(a)] $n_{ih}$ are non-negative integers,

\item[(b)] sums of entries in each row (from top to bottom) of $N$
equal the entries of $F-D$ (recall that rank $F$ is $t$),

\item[(c)] sum of entries in each column (from left to right) of
$N$ equal the entries of $E$ (recall that rank $E$ is $s$), and

\item[(d)] $n_{ih}= 0 \quad \mbox{if}\quad  m_{ih}= 0$.

\end{enumerate}
When $N$ has few non-zero entries, the conditions (a) to (c) will
be rather restrictive. When there is only one such non-trivial
matrix possible (say, as in the case of $M(T_1)$), then the
coefficient of $\Pi_{i,h} \beta_{ih}^{n_{ih}} $ in
$\Delta_{(D,E,F),(J,B)}$ is simply
$$
\Delta_{(D,E,F), N} = \Delta_{(D,E,F), (J,B)} \mid_{\beta_{ih}=1
\mbox{ \small{if} } n_{ih}\neq 0, \mbox{ \small{and} }
\beta_{ih}=0 \mbox{ \small{if} } n_{ih}=0}  \qquad \qquad (\ast)
$$
In some cases, condition (d) allows an inductive way of writing
$\Delta_{(D,E,F), M(T)}$ as a sum of determinants of the type
$(\ast)$. This is the case for our example, and we can compute the
following
{\small
$$
\Delta_{(D,E,F), M(T)} = \left| \begin{array}{ccccccccc}
X_{5,3} &0 &0 &0 &0 &0 &0 &Y_{5,2} &Y_{5,1} \\
0      &X_{5,3}&0 &0 &0 &0 &Y_{5,3} &0 &0 \\
0 &0 &X_{4,2} &0 &0 &Y_{4,3} &0 &Y_{4,2} &0 \\
0 &0 &0 &X_{3,1} &0 &Y_{3,3} &Y_{3,3} &0 &0 \\
0 &0 &0 &0 &X_{1,1} &0 &0 &0 &0 \\
0 &0 &0 &0 &0 &Y_{1,3} &0 &0 &0 \end{array} \right|
$$}
and the next is a sum of three determinants: {\small
\allowdisplaybreaks
$$
\begin{aligned}
& \Delta_{(D,E,F), M(T_1)} = \left| \begin{array}{ccccccccc}
X_{5,3} &0      &0      &0 &0 &0               &Y_{5,3} &0 &Y_{5,1} \\
0       &X_{5,3}&0      &0 &0 &Y_{5,3} &0               &Y_{5,2} &0 \\
0       &0      &X_{4,2} &0 &0 &0              &Y_{4,3} &Y_{4,2} &0 \\
0       &0     &0      &X_{3,1} &0 &Y_{3,3} &Y_{3,3} &0 &0 \\
0      &0      &0      &0      &X_{1,1} &0     &0 &0 &0 \\
0      &0      &0      &0      &0           &Y_{1,3} &0 &0  &0
\end{array} \right| \\
& -\left| \begin{array}{ccccccccc}
X_{5,3} &0      &0      &0 &0 &0               &Y_{5,3} &0 &Y_{5,1} \\
0       &X_{5,3}&0      &0 &0 &Y_{5,3} &0               &0 &0 \\
0       &0      &X_{4,2} &0 &0 &0              &0 &Y_{4,2} &0 \\
0       &0     &0      &X_{3,1} &0 &0 &Y_{3,3} &0 &0 \\
0      &0      &0      &0      &X_{1,1} &0     &0 &0 &0 \\
0      &0      &0      &0      &0           &Y_{1,3} &0 &0  &0
\end{array} \right| -
\left| \begin{array}{ccccccccc}
X_{5,3} &0      &0      &0 &0 &0               &Y_{5,3} &0 &Y_{5,1} \\
0       &X_{5,3}&0      &0 &0 &0 &0               &Y_{5,2} &0 \\
0       &0      &X_{4,2} &0 &0 &0              &Y_{4,3} &0 &0 \\
0       &0     &0      &X_{3,1} &0 &Y_{3,3} &0 &0 &0 \\
0      &0      &0      &0      &X_{1,1} &0     &0 &0 &0 \\
0      &0      &0      &0      &0           &Y_{1,3} &0 &0  &0
\end{array} \right|
\end{aligned}
$$}
Note that $F= (5,5,4,3,1,1)$ gives the heights of the blocks,
from top to bottom. The widths of the blocks (from left to right) are $D=(3,3,2,1,1)$ for
the $x$'s and $E=(3,3,2,1)$ for the $y$'s.

\smallskip

\subsection{Relationship Between $T$ and $\Delta_{(D, E, F), M(T)}$}

The polynomials $\Delta_{(D, E, F), M(T)}$ are coefficients  of
the polynomial $\Delta_{(D, E, F), (J, B)}$ of equation
(\ref{8.2}) with respect to the entries $\beta_{ih}$ of $B$. This
is a rather indirect definition, and we want to establish a more
direct connection between the polynomials $\Delta_{(D, E, F),
M(T)}$ and their LR tableaux.

\smallskip

The polynomial $\Delta_{(D, E, F), (J, B)}$ is the determinant of
the matrix $\widetilde Z_{(D, E, F)} = \widetilde{Z_{\ }}$,
described in equation (\ref{5.8}). The first $|D|$ columns of
$\widetilde{Z_{\ }}$ consist of the matrix $\widetilde X$, which
is partitioned into the blocks $\widetilde X_{jk} = \alpha_{jk}
X_{f_j, d_k}$. Since we have specialized the coefficients
$\alpha_{jk} = \delta_{jk}$, where the $\delta_{jk}$ are
Kronecker's deltas, only one block in each column of $\widetilde
X$ is non-zero. The condition that $D \subseteq F$ says that $d_k
\leq f_k$, so each block $X_{f_k, d_k}$ is taller than it is wide.
Therefore, each term in the expansion of $\widetilde{Z_{\ }}$ must
contain one entry from each column of $X_{f_k, d_k}$ as a factor.

\smallskip

Consider the terms which have  $\prod _{j = 1} ^{d_k} x_{jj}$ as
factor, for all $k$ where $1 \leq k \leq r$. Recall that $r$ is
the depth of the Young diagram $D$ (see (\ref{5.3})).  This will
be the determinant of a submatrix of $\widetilde{Z_{\ }}$, the
submatrix obtained by eliminating all the rows and columns of the
diagonal entries of each submatrix $\widetilde X_{jj}$. The
columns eliminated are exactly the columns of the submatrix
$\widetilde X$ of $\widetilde{Z_{\ }}$.  Therefore, the submatrix
$\widetilde{Z_{\ }}$ in which we are interested is actually a
submatrix of $\widetilde{Y_{\ }}$, which we will denote
$\widetilde{Y_o}$.  The submatrix $\widetilde{Y_o}$ has a block
structure, which is subordinate to the block structure of
$\widetilde{Y_{\ }}$, in the sense that each block of
$\widetilde{Y_o}$ is a submatrix of a block of $\widetilde{Y_{\
}}$.

More precisely, the block $(\widetilde{Y_o})_{j,k}$ consists of
the bottom $f_j - d_j$ rows of $\widetilde{Y}_{j,k} = \beta_{jk}
Y_{f_j, e_k} $. We will denote this submatrix of
$\widetilde{Y}_{jk}$ by
\begin{equation}\label{16.1}
\beta_{jk} Y_{(d_j, f_j], e_k}.
\end{equation}
\noindent Here, the notation $(d_j, f_j]$ is
intended to connote the open interval from $d_j$ to
$f_j$, and indicates that we are including the rows $d_j +1$ through $f_j$.

\smallskip

Summarizing, we may say that the coefficient of $\prod_{k = 1} ^r
\left ( \prod_{j = 1} ^{d_k} x_{jj}\right )$ in $\Delta_{(D, E,
F), (J, B)}$ is the determinant of the matrix $\widetilde{Y_o}$,
and the coefficient of $\prod_{k = 1}^r \left ( \prod_{j = 1}
^{d_k} x_{jj}\right )$ in $\Delta_{(D, E, F), M(T)}$ is the
coefficient of the monomial $M(T)$ of formula (\ref{15.1}) in the
determinant of $\widetilde{Y_o}$, considered as a function of the
$\beta_{ih}$. We will write
$$
\det \widetilde{Y_o} = \Delta_{(D, E, F), B, Y}.
$$
We are interested in the expansion of $\Delta_{(D, E, F), B, Y}$
with respect to the coefficients $\beta_{ih}$:
\begin{equation}\label{16.2}
\Delta_{(D, E, F), B, Y} = \sum_M \Delta_{(D, E, F), M, Y} \left(
\prod_{M=(m_{ih})} \beta_{ih} ^{m_{ih}} \right) ,
\end{equation}
and in particular, we define using (\ref{15.1}) and (\ref{16.2}),
\begin{equation}\label{16.3}
\delta_{T, Y} = \Delta_{(D, E, F), M(T), Y}.
\end{equation}

\smallskip

\noindent {\bf Example:} In our example \S 2.3, for tableaux $T$,
$T_1$, $T_2$ and $T_3$, the coefficient of $\prod_{k = 1}^r \left
( \prod_{j = 1} ^{d_k} x_{jj}\right )$ in $\Delta_{(D,E,F), (J,
B)}$ is given by {\small
$$
\det \widetilde{Y_o} = \Delta_{(D, E, F), B, Y} =
\left| \begin{array}{cccc}
\beta_{11}Y_{(3,5],3} &\beta_{12}Y_{(3,5],3} &\beta_{13}Y_{(3,5],2} &\beta_{14}Y_{(3,5],1} \\
\beta_{21}Y_{(3,5],3} &\beta_{22}Y_{(3,5],3} &\beta_{23}Y_{(3,5],2} &\beta_{24}Y_{(3,5],1} \\
\beta_{31}Y_{(2,4],3} &\beta_{32}Y_{(2,4],3} &\beta_{33}Y_{(2,4],2} &\beta_{34}Y_{(2,4],1} \\
\beta_{41}Y_{(1,3],3} &\beta_{42}Y_{(1,3],3} &\beta_{43}Y_{(1,3],2} &\beta_{44}Y_{(1,3],1} \\
\beta_{61}Y_{(0,1],3} &\beta_{62}Y_{(0,1],3} &\beta_{63}Y_{(0,1],2} &\beta_{64}Y_{(0,1],1} \end{array} \right|
$$ }
We note that
$$
F-D= (5-3, 5-3, 4-2, 3-1, 1-1, 1-0)= (2,2,2,2,0,1)
$$
which gives the heights of the blocks, from top to bottom. The
widths of the blocks are $E=(3,3,2,1)$, from left to right.
\smallskip

\noindent The coefficient of $\prod_{k = 1}^r \left ( \prod_{j =
1} ^{d_k} x_{jj}\right )$ in $\Delta_{(D,E,F), M(T)}$ is the same
as the coefficient of $M(T)$ in $\det \widetilde{Y_o}$ and is
given by {\small
$$
\Delta_{(D,E,F), M(T), Y} = \delta_{T,Y} = \left|
\begin{array}{cccc}
0 &0 &Y_{(3,5],2} &Y_{(3,5],1} \\
0 &Y_{(3,5],3} &0 &0 \\
Y_{(2,4],3} &0 &Y_{(2,4],2} &0 \\
Y_{(1,3],3} &Y_{(1,3],3} &0 &0 \\
Y_{(0,1],3} &0 &0 &0
\end{array} \right|
$$ }

\subsection{The Polynomials $\delta_{T,Y}$ (in the
Variables $y_{ac}$) Associated to a Tableau $T$}

For a tableau $T$, we have defined a monomial $M(T)$ in the
variables $\beta_{ih}$. We then considered the coefficient
$\Delta_{(D, E, F), M(T)}$ of $M(T)$ in $\Delta_{(D, E, F), (J,
B)}$. This is a polynomial in the $x_{ab}$'s and the $y_{ac}$'s.
By considering a standard coefficient in the $x$-variables of
$\Delta_{(D, E, F), M(T)}$, we have found a polynomial
$\Delta_{(D, E, F), M(T), Y} = \delta_{T, Y}$ in the variables
$y_{ac}$'s. This polynomial is associated to the LR tableau $T$,
but in a rather indirect way. We would like to relate it more
directly to $T$.

\smallskip

This is what we will do. We will associate to $T$ a monomial $e(T)$ in the
variables $y_{ac}$'s. Also, we will  define a term order on monomials in the
variables $y_{ac}$'s. Then we will show that, with respect to the given term
order, the monomial $e(T)$ is  the highest order term in the polynomial
$\delta_{T, Y}$.

\smallskip

The LR tableau $T$ is specified
by giving the entry of each box of $F^t - D^t$. To each box $\bf b$ of $T$,
we associate the variable $y_{a({\bf b}) c({\bf b})}$,
where $a({\bf b})$ is the row of
$F^t$ in which the box $\bf b$ lies, and $c ({\bf b})$ is the entry in $\bf
b$. To the tableau
$T$, we associate the product

\begin{equation}\label{17.1}
e(T) = \prod_{{\bf b} \in T } y_{a({\bf b}) c({\bf b})} =
\prod_{a,c} y_{ac}^{\ell(a,c)},
\end{equation}
where $\bf b$ runs over all boxes of $T$. Also define the monomial
\begin{equation}\label{new}
{\mathcal E} (T) = e(T) \prod_{k = 1}^r \left ( \prod_{j = 1}
^{d_k} x_{jj}\right )
\end{equation}
as the monomial which contains $e(T)$, in the expansion of
$\Delta_{(D, E, F), M(T)}$ (see \S 3.1, (2.2) and (3.1)).

\noindent {\bf Example:} Using the example in \S 2.3, we have
$$
e(T) = y_{53}y_{43}y_{52}y_{42}y_{32}y_{41}y_{31}y_{21}y_{11}
$$
using the following association by filling up the skew diagram
$F^t - D^t$ with $y$'s: {\small
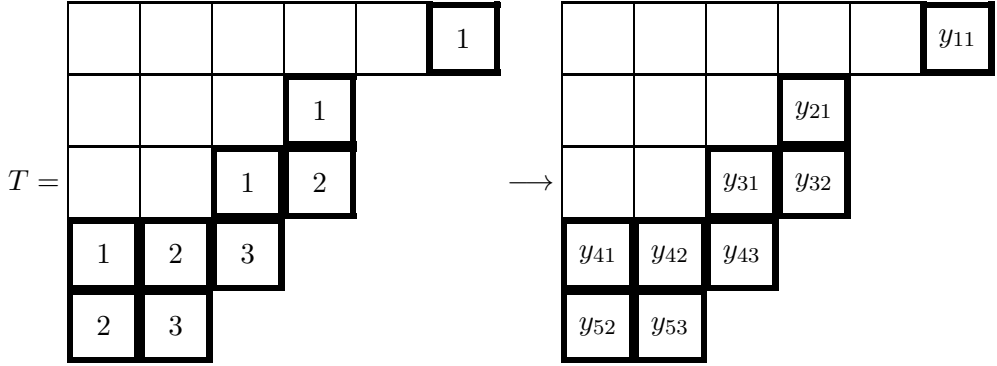
$$
T =\tableau[m]{ & & & & &\tf{1} \\ & & &\tf{1} \\ & &\tf{1} &\tf{2} \\
\tf{1} &\tf{2} &\tf{3} \\ \tf{2} &\tf{3} }  \longrightarrow
\tableau[m]{ & & & & &\tf{y_{11}} \\ & & &\tf{y_{21}} \\ & &\tf{y_{31}} &\tf{y_{32}}\\ \tf{y_{41}} &\tf{y_{42}} &\tf{y_{43}} \\ \tf{y_{52}} &\tf{y_{53}} }
$$
Compare with $e(T_1) = y_{53}^2 y_{42}^2 y_{32}y_{41}y_{31}y_{21}y_{11}$:
$$
T_1 =\tableau[m]{ & & & & &\tf{1} \\ & & &\tf{1} \\ & &\tf{1} &\tf{2} \\ \tf{1} &\tf{2} &\tf{2} \\ \tf{3} &\tf{3} }
\longrightarrow
\tableau[m]{ & & & & &\tf{y_{11}} \\ & & &\tf{y_{21}} \\ & &\tf{y_{31}} &\tf{y_{32}}\\ \tf{y_{41}} &\tf{y_{42}} &\tf{y_{42}} \\ \tf{y_{53}} &\tf{y_{53}} }
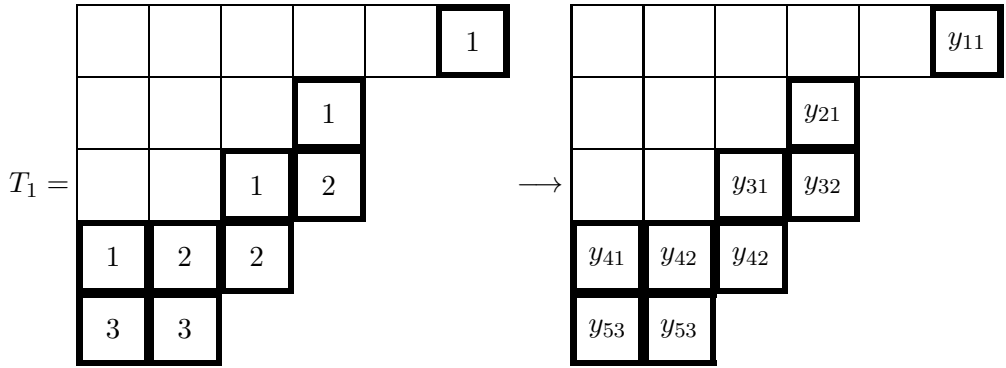
$$
}

\medskip

\begin{lemma}\label{lem18}
The LR tableau $T$ is determined by the monomial $e(T)$ together
with the triple $(D, E, F)$. Further, the map $T \mapsto {\mathcal
E} (T)$ is one to one.
\end{lemma}

\smallskip

\noindent {\bf Proof:} Indeed, we can see
from the definition (3.5) of $e(T)$ that it may also be written
$$
e(T) = \prod_{a,c} y_{ac}^{\ell(a,c)},
$$
where $\ell(a, c)$ is the number of entries equal to $c$ in row
$a$. In other words, $e(T)$ records the number of boxes of each
possible entry in each row. Since the entries are always arranged
in weakly increasing order, and there are no gaps between
consecutive boxes, this means that $e(T)$ determines each row of
$T$, up to translation in the horizontal direction. But, we also
know that $T$ occupies the boxes of $F^t - D^t$, so that its
location is also fixed. This proves the first statement of the
lemma.

For the second statement, note that the factor $\prod_{k = 1}^r
\left ( \prod_{j = 1} ^{d_k} x_{jj}\right )$ in the monomial
${\mathcal E}(T)$ (see (\ref{new})) encodes the Young diagram
$D^t$. Now, the banal tableau of shape $E^t$ and the filling of
$F^t -D^t$ can be recovered from $e(T)$. The LR tableau $T$ is
thus  effectively determined by ${\mathcal E}(T)$. Hence, the
second statement. $\qquad \square$

\subsection{Proof of the Theorem}

\subsubsection{Monomial Ordering in the $y$ Variables}
We now fix a monomial ordering on monomials
in the variables $y_{ac}$'s. First, we specify that
$$
y_{ac} > y_{a'c'}
$$
if and only if $c' > c$, or $c = c'$ and $a' > a$. Note that this
is a total ordering on the $y_{ac}$'s. It is essentially lexicographic order,
except that we order on the second index first. Thus,
$$
y_{11} > y_{21} > \ldots > y_{n1} >  y_{12} > y_{22} > \ldots > y_{n\ell}.
$$

We extend this to an ordering between all monomials on the
$y_{ac}$'s by graded lexicographic order (see [CLO]).  That is, if
$M_1$ and $M_2$ are two monomials of unequal degree, then the one
with larger degree is larger. If they have the same degree, we
write them as products

$$
M_1 = y_{a_1 c_1} y_{a_2 c_2} y_{a_3 c_3} \cdots \cdots ,
$$
and
$$
M_2 = y_{a'_1 c'_1} y_{a'_2 c'_2} y_{a'_3 c'_3} \cdots \cdots,
$$
where the factors of each monomial are written in decreasing
order. Suppose that $k$ is the first index such that $y_{a_k c_k}
\neq y_{a'_k c'_k}$. Then $M_1 > M_2$ if and only if $y_{a_k c_k}
> y_{a'_k c'_k}$.

\medskip

\begin{lemma}\label{lem20}
Given an LR tableau $T$, associate to it the polynomial \linebreak
$\Delta_{(D, E, F), M(T), Y }= \delta_{T, Y}$ in the variables
$y_{ac}$'s, as described in formula (\ref{16.3}). Then the largest
monomial, with respect to the monomial ordering described above,
in the expansion of $\delta_{T, Y}$ is the monomial $e(T)$ defined
in formula (\ref{17.1}).
\end{lemma}

\smallskip

The proof of Lemma \ref{lem20} is by induction, and will be
organized in several parts. It will occupy \S 3.4.2 to \S 3.4.6.

\subsubsection{The Nature of $\widetilde{Y_o}$}

To begin the proof of Lemma \ref{lem20}, we consider
the nature of the matrix $\widetilde{Y_o}$ of \S 3.2, and the role
of the monomial $e(T)$ in determining $\delta_{T,Y}$.

\smallskip

The matrix $\widetilde{Y_o}$ is partitioned into blocks. The blocks
of columns correspond to columns of $E^t$; their sizes are $e_k$,
the parts of the partition $E$. The blocks of rows correspond to
the columns of $F^t - D^t$; their sizes  are $f_j - d_j$, the
differences between the parts of $F$ and the parts of $D$. The
matrices which occupy the blocks are the matrices
$(\widetilde{Y_o})_{j,k} = \beta_{jk} Y_{(d_j, f_j], e_k}$ of
equation (\ref{16.1}). In the discussion below, we will have to be referring
both to individual rows or columns of the matrix $\widetilde{Y_o}$,
and to rows or columns of blocks of $\widetilde{Y_o}$. To alleviate
somewhat the awkwardness of distinguishing between the two, we
will refer to the  blocks of rows as \lq \lq {\it superrows}" and
to the blocks of columns  as \lq \lq {\it supercolumns}". Thus,
the matrix $(\widetilde{Y_o})_{j,k}$, considered as a submatrix of
$\widetilde{Y_o}$, is the intersection of the $j$-th superrow and
the $k$-th supercolumn.

\smallskip

Consider the usual alternating sum expansion of
$\det \widetilde{Y_o}$, with each term being a product of entries of
$\widetilde{Y_o}$, one from each row and column. We see from the form
(\ref{16.1}) of the blocks of $\widetilde{Y_o}$, that a given term
in this expansion will have a factor of $\beta_{ih}$ for each
entry chosen from the $(i, h)$-block. Thus, for a given monomial
$\prod_{i, h} \beta_{ih}^{m_{ih}}$, the exponent $m_{ih}$
determines the number of entries in the $(i, h)$-block that a term
with this $\beta$-monomial will have.  The total number of terms
chosen from a given superrow of $\widetilde{Y_o}$ will have to be
the number of rows in that superrow; and similarly for columns.
From this observation, it follows that, for a monomial $\prod_{i,
h} \beta_{ih}^{m_{ih}}$ to have a non-zero coefficient in the
expansion of $\det \widetilde{Y_o}$, the conditions for row sums and column sums, respectively,
\begin{equation}\label{21.1}
\sum_{h = 1}^s  m_{ih} = f_i - d_i \quad \quad {\rm and} \quad \quad
\sum_{i = 1}^t  m_{ih} = e_h
\end{equation}
must hold.

\smallskip

Also, we may observe that, if some $\beta_{ih}$ occurs does not occur in
a monomial
$M(T)$ (that is, it occurs with exponent zero), then any term in the coefficient
of $M(T)$ will have to completely avoid the block $(\widetilde{Y_o})_{i,h}$. This
means that, as far as the coefficient of $M(T)$ is concerned, we will get the
same result if we replace $(\widetilde{Y_o})_{i,h}$ with the zero matrix.

\smallskip

\noindent {\bf Example:} Reverting to our same tableau $T$ as in
\S 2.3, we look at the exponent matrix on the left. Note that
$m_{5j}=0$ because $f_5 -d_5=1-1=0$.  Because in the determinant
expansion, we need to select one and only one element from each
column as well as each row, we see that the row sums must be
$F-D=(2,2,2,2,0,1)$ and the column sums must be $E=(3,3,2,1)$ as in (3.7):
{\small
$$
\left[ \begin{array}{cccc}
0 &0 &1 &1 \\
0 &2 &0 &0 \\
1 &0 &1 &0 \\
1 &1 &0 &0 \\
0 &0 &0 &0 \\
1 &0 &0 &0
\end{array} \right]
\quad \leftrightarrow \quad
\det \widetilde{Y_o} =\left| \begin{array}{cccc}
\beta_{11}Y_{(3,5],3} &\beta_{12}Y_{(3,5],3} &\beta_{13}Y_{(3,5],2} &\beta_{14}Y_{(3,5],1} \\
\beta_{21}Y_{(3,5],3} &\beta_{22}Y_{(3,5],3} &\beta_{23}Y_{(3,5],2} &\beta_{24}Y_{(3,5],1} \\
\beta_{31}Y_{(2,4],3} &\beta_{32}Y_{(2,4],3} &\beta_{33}Y_{(2,4],2} &\beta_{34}Y_{(2,4],1} \\
\beta_{41}Y_{(1,3],3} &\beta_{42}Y_{(1,3],3} &\beta_{43}Y_{(1,3],2} &\beta_{44}Y_{(1,3],1} \\
\beta_{61}Y_{(0,1],3} &\beta_{62}Y_{(0,1],3} &\beta_{63}Y_{(0,1],2} &\beta_{64}Y_{(0,1],1} \end{array} \right|
$$}

\subsubsection{Further Constraints on $M(T)$}

For a collection of exponents $m_{ih}$ to be the exponents of a
monomial $M(T)$ coming from an LR tableau $T$, then the $m_{ih}$
must satisfy further constraints besides the equations
(\ref{21.1}). Recall from \S 3.1 that the exponent $m_{ih}$ of the
monomial $M(T)$ is the number of boxes from the $h$-th column of
$BE^t$ which get put into the $i$-th column of $F^t$ by the
inverse of the standard peeling map. From the description of the
standard peeling in \S 2.3, we see that in the filling of $F^t$,
the $m$-th element of column $i$ of $BE^t$ lies to the northeast,
in fact, strictly right and weakly above, the $m$-th element of
column $i+1$. (Note that the $m$-th element of column $i$ of
$BE^t$ is the number $m$, since $BE^t$ is banal.) It follows that
the number of elements of column $i$ of $BE^t$ which are assigned
to columns $k+1$ or greater of $F^t$ is always at least as large
as the number of elements from column $i+1$ of $BE^t$ which are
assigned to columns $k$ or greater of $F^t$. We can state these
conditions in terms of the numbers $m_{ih}$:
$$
\sum_{j = k+1}^t m_{ji}  \geq \sum_{j = k}^t m_{j,(i+1)}.
$$
In particular, if $i(h)$ is the column of $F^t$ which gets the
first entry of column $h$ of $BE^t$,  that is, column $i(h)$ is
the rightmost column of  $F^t$ to receive elements from column
$h$, then $i(h)$ is strictly decreasing with $h$.

\subsubsection{Acquiring Factors $y_{a1}$'s in the Leading Term in the Coefficient $\delta_{T, Y}$ of $M(T)$}

Now fix a tableau $T$, and consider the leading term in the
coefficient $\delta_{T, Y}$ of $M(T)$ in the expansion (\ref{16.3}) of
$\det \widetilde{Y_o}$. Since the most important variables in
determining the ordering of the term of $\delta_{T, Y}$ are the
$y_{a1}$'s, we will focus first on understanding how terms in
$\delta_{T,Y}$ can acquire these variables as factors.

As remarked at the end of \S 3.4.2, for any $\beta_{ih}$ which
occurs in $M(T)$ with exponent zero, we may replace the block
$(\widetilde{Y_o})_{i,h}$ with zeroes without changing $\delta_{T,
Y}$. We will do this.  Call the resulting matrix
$\widetilde{Y_T}$. Since the first element from the $(i+1)$-th
column of $BE^t$ is always placed strictly to the left of the
first element of the $i$-th column, we see that the matrix
$\widetilde{Y_T}$ will have a block triangular form. More
precisely, in the $h$-th supercolumn of $\widetilde{Y_T}$, the
last non-zero block of $\widetilde{Y_T}$ will be in a strictly
lower superrow than the last non-zero block of the $(h+1)$-th
supercolumn. In the notation of \S 3.4.3, the last block of the
$h$-th supercolumn is $(\widetilde{Y_T})_{i(h),h }$. The
triangularity of $\widetilde{Y_T}$ amounts to the fact that
$i(h+1) < i(h)$. We will call the last non-zero block in a given
supercolumn of $\widetilde{Y_T}$ the \lq \lq final block" of that
supercolumn.

\smallskip

 What does this imply about the occurrence of the variables
$y_{a1}$'s in the terms in the expansion of $\det \widetilde{Y_T}$?
The variables $y_{a1}$'s appear in the first (leftmost) column of
a supercolumn. This makes $s$ columns all together which contain
entries $y_{a1}$'s. We will call these columns the \lq \lq {\it
ones columns}". Any term in the standard expansion of $\det
\widetilde{Y_T}$ will select one factor from each ones column, so
the sum of all the exponents of the variables $y_{a1}$'s in $M(T)$
must be $s$.

\smallskip

For the leading term, we want to make the indices $a$ in the
factors $y_{a1}$'s as small as possible. We will refer to the ones
column of the $h$-th supercolumn as the $h$-th ones column. Let
$a_h$ be the smallest index among the indices $a$ making up the
$y_{a1}$'s belonging to the $h$-th ones column. Then the largest
conceivable product of elements chosen one from each ones columns
is
\begin{equation}\label{23.1}
e_1(T) = \prod_{h = 1}^s y_{a_h 1}.
\end{equation}

\noindent {\bf Example:} We recall the matrix $\widetilde{Y_T}$, in
the case for our $T$ as in \S 2.3.  Notice that the final blocks
$Y_{(0,1],3}$, $Y_{(1,3],3}$, $Y_{(2,4],2}$, $Y_{(3,5],1}$ slopes
upwards from left to right. Also since $E=(3,3,2,1)$, i.e., $s=4$,
we have four supercolumns with blocks containing $y_{a1}$'s.
$$
\det \widetilde{Y_T} = \left|
\begin{array}{cccc}
0 &0 &Y_{(3,5],2} &Y_{(3,5],1} \\
0 &Y_{(3,5],3} &0 &0 \\
Y_{(2,4],3} &0 &Y_{(2,4],2} &0 \\
Y_{(1,3],3} &Y_{(1,3],3} &0 &0 \\
Y_{(0,1],3} &0 &0 &0
\end{array} \right|
$$

\smallskip

\noindent To form $e_1(T)$, pick $y_{11}$ from the final block $Y_{(0,1],3}$ in the
first supercolumn, pick $y_{21}$ from the final block
$Y_{(1,3],3}$ in the second supercolumn, pick $y_{31}$ from the
final block $Y_{(2,4],2}$ in the third supercolumn and finally
pick $y_{41}$ from the final block $Y_{(3,5],1}$ in the last
supercolumn.  Thus $e_1(T) = y_{11}y_{21}y_{31}y_{41}$.

\smallskip

In fact we can realize the product $e_1 (T)$ as a factor in
elements in the term expansion of $\det \widetilde{Y_T}$, and it is
not hard to see explicitly how to do this.  Indeed, we can
explicitly locate one occurrence of $y_{a_h 1}$ in the $h$-th ones
column. It is in the top row of the final block of the
supercolumn.  This follows from the formula (\ref{16.1}) defining
the blocks of $\widetilde{Y_o}$ (and thus defining the non-zero
blocks of $\det \widetilde{Y_T}$). It is clear from this formula that if the
$(i, h)$ block of $\widetilde{Y_T}$ is non-zero, its ones column
contains entries $y_{a1}$ for $a$ running from $d_i +1$ to $f_i$.
Since the $d_i$ decrease as $i$ increase, it follows that the
smallest possible first entry for the $y_{a1}$ in a given ones
column will occur in the final block of the corresponding
supercolumn. (It may occur also in blocks above the last one.) Let
us call this occurrence of $y_{a_h 1}$ the \lq \lq {\it reference
entry}" of the $h$-th ones column.

\smallskip

Because of the triangular character of $\widetilde{Y_T}$, the final
blocks of different supercolumns occur in different superrows. We
conclude that the reference entries of different ones columns
occur in different rows, and therefore, the product $M_1(T)$ can
occur in the terms of the standard expansion of
$\det \widetilde{Y_T}$. To prove Lemma \ref{lem20}, we will investigate which terms
of the expansion of $\det \widetilde{Y_T}$ will have $e_1(T)$ as a
factor.

\subsubsection{Factoring $\det \widetilde{Y_T}$ by $e_1 (T)$}

Several ones columns may have the same variable $y_{a 1}$ as their
reference entry. For a given $a$, consider the set of $h$ such
that $a_h = a$, that is, $y_{a 1} = y_{a_h 1}$ is the reference
entry of the $h$-th ones column. Suppose that there are $m'_{a1}$
such ones columns. The supercolumns containing these ones columns
correspond to columns of $BE^t$ whose first entry gets put in the
$a$-th row of $F^t$. From the definition of the standard peeling
and condition (LR1) defining an LR tableau, we see that this will
be a consecutive set of columns of $BE^t$, and that the ones
entries of these columns will fill consecutive boxes in the $a$-th
row of $F^t$. These boxes will be a consecutive set of $m'_{a1}$
boxes in the $a$-th row of  $F^t$; in fact, they will be the
leftmost $m'_{a1}$ boxes in the $a$-th row of $F^t - D^t$. The
next column of $F^t - D^t$ to the left of the leftmost of the
column of this set will start in a row below the $a$-th row of
$F^t$.  As an illustration, look at the entries $y_{11}$,
$y_{21}$, $y_{31}$ and $y_{41}$ in the fillings of $F^t - D^t$ in
the Example of \S 3.3.

\smallskip

We will call the first element of the first column of a given block, the
$(1,1)$ entry of the block. If the $(1, 1)$ entry of a given block is
$y_{a1}$, this means that the column of $F^t - D^t$ corresponding to the
superrow of the block starts in the $a$-th row of $F^t$.

\smallskip

We may conclude from this discussion that the following statements hold.
\begin{enumerate}
\item[(i)] The ones columns which have a given variable $y_{a1}$ as their reference entry form a consecutive set of ones columns, say from $h_0(a)$ to $h_1(a)$.  (Hence,
$h_1 (a) - h_0 (a) + 1 = m'_{a1}$.)

\item[(ii)] If the supercolumn $h$, with  $h_0(a) \leq h < h_1 (a)$, has final block in the $i$-th superrow, then supercolumn $h+1$ has final block in superrow $(i-1)$.

\item[(iii)] For column $h$, with $h_0(a) \leq h \leq h_1 (a)$, the number of $(1,1)$  entries which are equal to $y_{a 1}$ is equal to $(h_1 (a) -  h +1)$. These are the lowest $(h_1(a) - h +1)$ $(1,1)$ entries on the $h$-th superrow.
\end{enumerate}

\smallskip

On the basis of  facts (i), (ii), (iii) above, we may
make the following conclusions.

\begin{enumerate}

\item[(a)] In a supercolumn $h$ with reference element $y_{a1}$,
we may add suitable multiples of the first row of the final block
to the first rows of the $(h_1(a) - h)$ blocks above the final
block, and make these rows all equal to zero. (This procedure
unfolds in the most orderly fashion if one does the row operations
from left to right, that is, beginning with $h = h_o(a)$ and
progressing consecutively to the larger $h$'s. Then the first rows
of the blocks in the supercolumns with $y_{a1}$ as reference
element lying  to the left of the first row of the final block of
the $h$-row will have been eliminated before one does row
operations in the $h$-th supercolumn, so these row operations will
leave the earlier supercolumns unchanged. Alternatively, one could
work backwards, moving to the right from the $h_1(a)$-th
supercolumn. The row operations in supercolumn $h$ will then alter
the supercolumns to the left, but whatever effects it has can be
cleaned up when it comes time to do row operations in those
supercolumns.) In particular, these row operations, which will not
affect $\det \widetilde{Y_T}$, will leave only one entry equal to
$y_{a1}$ in any of these columns, namely, the reference element
itself.

\item[(b)] Moreover, all these row operations affect only elements which share a
row with one of the reference elements. Hence, they do not affect the
cofactor of the reference elements.
\end{enumerate}

\smallskip

After we have finished with the operations of statement (a) above for all values of $a$ (as the first index of $y_{a1}$), we are left with a matrix such that the only entry in any ones column equal to the reference entry is the
reference entry itself. This leads us to the following statement.

\medskip

\begin{lemma}\label{sublem26}
With the monomial $e_1(T)$ as in formula (\ref{23.1}), we have
$$
\det \widetilde{Y_T} = e_1 (T) \Gamma + R,
$$
where
\begin{enumerate}
\item[(a)] $\Gamma$ is the cofactor in $\det \widetilde{Y_T}$ of the reference elements -- the determinant of the submatrix  of the rows and columns {\bf not} containing the reference elements; and
\item[(b)] all terms of $R$ are dominated by any term of $e_1(T) \Gamma$.
\end{enumerate}
\noindent In particular, the leading term of $\det \widetilde{Y_T}$
is $e_1(T)$ times the leading term of $\Gamma$.
\end{lemma}

\subsubsection{The Inductive Step}

Let's compare the monomial $e_1(T)$ of formula (\ref{23.1})
with the monomial $e(T)$ of formula (\ref{17.1}). Observe
that the first indices of the variables $y_{ac}$ which occur in a
given superrow label the rows of the column of $F^t - D^t$
corresponding to that superrow. These first indices are the same
for every block in the superrow. In particular, the index $a_h$ of
the reference element in the $h$-th ones column is the smallest
row in the $i(h)$-th column of $T= F^t - D^t$. This in turn is where
the 1 from the $h$-th column of $E^t$ gets put by the inverse of
the standard stripping of $T$. This means that the factor of the
monomial $e(T)$ corresponding to the first box of the $i(h)$-th
column of $F^t - D^t$ is just $y_{a_h 1}$. We conclude that
$e_1(T)$ divides $e(T)$. On the other hand, we know that $e_1 (T)$
has degree $s$ equal to the number of columns of $E^t$. From our
description of the standard stripping, we know that this is the
same as the number of boxes of $T$ containing a 1. Therefore, the
monomial $e_1(T)$ contains exactly all the factors $y_{a1}$
dividing $e(T)$: the quotient $e(T)/e_1(T)$ is devoid of the
variables $y_{a1}$ -- it contains only variables $y_{ac}$ with $c
\geq 2$. We may also express this as follows: the monomial
$e_1(T)$ is the factor of the monomial $e(T)$ attributable to the
boxes of $T$ filled with a 1, and the quotient $e(T)/e_1 (T)$ is
the factor attributable to the boxes filled with the numbers 2 or
larger.

\smallskip

Now consider the matrix whose determinant is the factor $\Gamma$
in the Lemma \ref{sublem26}. It is the matrix obtained from
$\widetilde{Y_T}$ be removing all the first columns of each
supercolumn, and all the first rows of the superrows containing
final blocks. We will call this matrix $(\widetilde{Y_T})_{/1}$. It
inherits a partitioned structure from $\widetilde{Y_T}$. The blocks
of $(\widetilde{Y_T})_{/1}$ are just the blocks of $\widetilde{Y_T}$
with the first column removed, and possibly also with the first
row removed, if the block shares a row with a final block.

\smallskip

Let $T_{/1}$ be the subtableau of $T$ obtained by
omitting all the cells of $T$ which contain a 1. Checking the
criteria  of \S 2.3 for $LR$ tableau, we see that, except for the
fact that $T_{/1}$ contains only entries $2$ through $\ell_0$, it
is an $LR$ tableau. More precisely, if we subtract 1 from each
entry of $T_{/1}$, it becomes an $LR$ tableau.

\smallskip

Now, by comparing the tableau $T_{/1}$, the matrix
$(\widetilde{Y_T})_{/1}$, and the monomial $e_{/1}(T) = e(T)/e_1 (t)$, we see
that they are parallel to the original triple of tableau $T$, matrix
$\widetilde{Y_T}$ and monomial $e(T)$. The matrix $(\widetilde{Y_T})_{/1}$  and the
monomial $e_{/1}(T)$ are related to the tableau $T_{/1}$ in the same way
that the matrix $\widetilde{Y_T}$ and the monomial $e(T)$ are related to the
original LR tableau $T$. Further all these relations are consistent with the
use only of the numbers 2 through $\ell_0$ to fill the boxes of $T$, rather
than starting with the number 1.

From these observations, we conclude that repetition of the
reasoning of \S 3.4.2 to \S 3.4.5 for $2, 3, \ldots, \ell_0$ proves
Lemma \ref{lem20}.

\subsubsection{Proof of Main Result}

Lemma \ref{lem20} makes it fairly easy to prove the main result of
this paper:

\smallskip

\noindent \bf Theorem: \rm {\it As the coefficient matrices $A$
and $B$ of formulas (\ref{5.4}) vary through all possible
constants, the polynomials $\Delta_{(D,E, F), (A, B)} (X, Y)$ of
formula (\ref{5.9}) span the tensor product algebra $TA_{n, k,
\ell}$. More precisely the polynomials $\Delta_{(D, E, F), M(T)}$
of formula (\ref{8.2}), with the monomial $M(T)$ given by
(\ref{15.1}), form a basis for $TA_{n, k, \ell}$.}

\smallskip

\noindent {\bf Proof:} We know from Lemma \ref{lem6} and the
formula (\ref{8.2}) that the polynomials $\Delta_{(D, E, F),
M(T)}$ are $GL_n$-highest weight vectors with weight $F^t$ in the
tensor product of the representations $\rho_n^{D^t}$ and $\rho_n
^{E^t}$ of $GL_n$. From Lemma \ref{lem20}, we know that the
leading term in $\Delta_{(D, E, F), M(T)}$ is the monomial $e(T)$
of formula (\ref{17.1}). Lemma \ref{lem18} tells us that $e(T)$
together with the diagrams $D$, $E$, and $F$ determine $T$. It
follows that the collection of polynomials $\Delta_{(D, E, F),
M(T)}$ corresponding to some collection of tableaux $T$, is
linearly independent. Since we know from [Ful] that LR tableaux
attached to the diagrams $D$, $E$ and $F$ count the multiplicity
of $\rho_n ^{F^t}$ in $\rho_n^{D^t} \otimes \rho_n ^{E^t}$, we see
that the $\Delta_{(D, E, F), M(T)}$ will be a basis for the
$\psi^{F^t} \times \psi^{D^t} \times \psi^{E^t}$ weight space of
$TA_{n, k, \ell}$. Letting the diagrams $F$, $E$ and $D$ vary, the
theorem follows. $\qquad \square$

\section{Example: ${\bf SL_4}$ Tensor Product Algebra}

\setlength{\unitlength}{2.5mm}

In this section, we will illustrate our results using the example
of the $SL_4$ tensor product algebra. Several people have worked
on this particular case [BZ], [Gro], [How], [Van].

Berenstein and Zelevinsky [BZ] approached this in terms of triple
multiplicities.  The multiplicities are described as the number of
integral points in certain convex sets.  Let ${\mathfrak
R}(SL_n/U_n)$ denote the algebra of regular functions on the
natural torus bundle over the flag manifold for $SL_n$. The
approach identifies the highest weight vectors in $(\mathfrak R
(SL_n/U_n) \otimes \mathfrak R (SL_n/U_n) \otimes \mathfrak R
(SL_n/U_n))^{SL_n}$ using Berenstein-Zelevinsky configurations
(we'll simply call them BZ diagrams), which are arrays of hexagons
and triangles, with integers at each vertex so that the sums on
opposite sides of any hexagon is the same. The generators for the
tensor product algebra correspond to the primitive BZ diagrams;
they comprise of 0's and 1's. We have drawn them in the following
diagram, with a ${\Large \bullet}$ referring to an integer 1 at
that vertex, and 0's at other vertices.

One important information is the $(\widehat{A_3}^+)^3$ grading
corresponding to the triple of diagrams $(D, E, F)$.  There are
several conventions to read this. We adopt the following: In the
case of $SL_4$, one will find numbers $(x_{ij}, y_{ij}, z_{ij})$
associated to 6 triangles in the BZ diagram:
$$
\begin{array}{ccccccccccc}
&&&&&x_{11}&&&&& \\
&&&&y_{11}&&z_{11}&&&& \\
&&&x_{12}&&&&x_{13}&&& \\
&&y_{12}&&z_{12}&&y_{13}&&z_{13}&& \\
&x_{21}&&&&x_{22}&&&&x_{23}& \\
y_{21}&&z_{21}&&y_{22}&&z_{22}&&y_{23}&&z_{23}
\end{array}
$$
The $(\widehat A_3^+)^3$ grading corresponding to this BZ diagram is $(D,E,F)$ where
$$
\begin{array}{ccc}
D^t&= &( x_{11}+ y_{11}+ x_{12}+ y_{12}+ x_{21}+ y_{21}, x_{11}+ y_{11} + x_{12}+ y_{12}, x_{11}+ y_{11} ) \\
E^t&= &( y_{21}+ z_{21}+ y_{22}+ z_{22}+ y_{23}+ z_{23}, y_{21}+ z_{21} + y_{22}+ z_{22}, y_{21}+ z_{21} ) \\
F^t&= &( x_{11}+ z_{11}+ x_{13}+ z_{13}+ x_{23}+ z_{23}, x_{11}+ z_{11} + x_{13}+ z_{13}, x_{11}+ z_{11} )
\end{array}
$$

\smallskip

\noindent {\bf Example:} Take the following primitive BZ diagram,
where each bold ``dot'' in the triangles represent ``1'', and zero
otherwise. This BZ diagram corresponds to $D^t=(2,1,1)$,
$E^t=(1,1,0)$ and $F^t=(1,1,0)$:

$$
\begin{picture}(10,7.5)
\thicklines
\put (0,0){\line(1,0){10}}
\put (0,0){\line(3,5){5}}
\put (4,0){\line(3,5){3}}
\put (8,0){\line(3,5){1}}
\put (2,0){\line(-3,5){1}}
\put (6,0){\line(-3,5){3}}
\put (10,0){\line(-3,5){5}}
\put (2.1,3.464){\line(1,0){5.8}}
\put (4.1,6.928){\line(1,0){1.8}}
\thinlines
\put (1,1.732){\line(1,0){2}}
\put (4,0){\line(-3,5){1}}
\put (3,1.732){\line(3,5){2}}
\put (7,5.196){\line(-1,0){2}}
\put (4,6.928){\line(3,-5){1}}
\put (1,1.732){\circle*{1}}
\put (4,3.464){\circle*{1}}
\put (4,6.928){\circle*{1}}
\put (4,0){\circle*{1}}
\put (7,5.196){\circle*{1}}
\end{picture}
$$

\smallskip

In the following table, we provide the generators of the tensor
product algebra for $SL_4$ using the primitive BZ diagrams. We
also provide in the last two columns, key information such as the
lead monomial $e(T)$ in $\delta_{T,Y}$ (see (3.5)) as well as the
monomial ${\mathcal E}(T)= e(T)\, \Pi_{k=1}^r \left(
\Pi_{j=1}^{d_k} x_{jj} \right)$ containing the monomial $e(T)$ in
the polynomial $\Delta_{(D, E, F), M(T)}$ of formula (\ref{8.2}).
Of course, this monomial ${\mathcal E}(T)$ determines $D^t$, $E^t$
and $F^t$ (see proof of Lemma 3.1). The exponents of $\Pi_{k=1}^r
\left( \Pi_{j=1}^{d_k} x_{jj} \right)$ determines $D^t$, while
$e(T)$ encodes the banal tableau $BE^t$ of shape $E^t$ and the
filling of the skew tableau $F^t - D^t$ (see Lemma 3.1).

\smallskip

\begin{tabular}{|c|c|c|c|c|c|c|c|} \hline\hline
\multicolumn{8}{|c|}{Generators of $SL_4$ Tensor Product Algebra}
\\ \hline\hline
{No.} &{$\ $ BZ Diagrams $\ $} &\ $D^t\ $  &\ $E^t\ $  &\ $F^t\ $
&$\Delta_{(D,E,F), M(T)}$ &$e(T)$ &${\mathcal E}(T)$ \\
\hline\hline &&&&&&& \\

${\begin{array}{c} 1 \\ \quad \\ \quad \\ \end{array}}$
& {\begin{picture}(10,7.5)
\thicklines
\put (0,0){\line(1,0){10}}
\put (0,0){\line(3,5){5}}
\put (4,0){\line(3,5){3}}
\put (8,0){\line(3,5){1}}
\put (2,0){\line(-3,5){1}}
\put (6,0){\line(-3,5){3}}
\put (10,0){\line(-3,5){5}}
\put (2.1,3.464){\line(1,0){5.8}}
\put (4.1,6.928){\line(1,0){1.8}}
\put (5,8.6){\circle*{1}}
\end{picture} }
&{$\tableau[b,s]{ \\  \\  \\ }$} &$\emptyset$ &{$\tableau[b,s]{\\ \\  \\
}$} &${\begin{array}{c} |X_{3,3}| \\ \quad \\ \quad \\
\end{array}} $
&${\begin{array}{c} 1 \\ \quad \\ \quad \\
\end{array}}$
&${\begin{array}{c} x_{11}x_{22}x_{33} e(T) \\ \quad \\ \quad \\
\end{array}}$
\\ \hline\hline  &&&&&&& \\

${\begin{array}{c} 2 \\ \quad \\ \quad \\ \end{array}}$
& {\begin{picture}(10,7.5)
\thicklines
\put (0,0){\line(1,0){10}}
\put (0,0){\line(3,5){5}}
\put (4,0){\line(3,5){3}}
\put (8,0){\line(3,5){1}}
\put (2,0){\line(-3,5){1}}
\put (6,0){\line(-3,5){3}}
\put (10,0){\line(-3,5){5}}
\put (2.1,3.464){\line(1,0){5.8}}
\put (4.1,6.928){\line(1,0){1.8}}
\put (0,0){\circle*{1}}
\end{picture} }
&{$\tableau[b,s]{\\ }$} &{$\tableau[b,s]{\\ \\ \\ }$}
&{$\tableau[b,s]{ \\ 1\\ 2\\ 3 }$} &${\begin{array}{c}
|X_{4,1}Y_{4,3}| \\ \quad \\ \quad \\ \end{array}}$
&${\begin{array}{c} y_{21}y_{32}y_{43} \\ \quad \\ \quad \\
\end{array}}$
&${\begin{array}{c} x_{11}e(T) \\ \quad \\ \quad \\
\end{array}}$
\\ \hline\hline  &&&&&&& \\

${\begin{array}{c} 3 \\ \quad \\ \quad \\ \end{array}}$
&{\begin{picture}(10,7.5) \thicklines \put (0,0){\line(1,0){10}}
\put (0,0){\line(3,5){5}} \put (4,0){\line(3,5){3}} \put
(8,0){\line(3,5){1}} \put (2,0){\line(-3,5){1}} \put
(6,0){\line(-3,5){3}} \put (10,0){\line(-3,5){5}} \put
(2.1,3.464){\line(1,0){5.8}} \put (4.1,6.928){\line(1,0){1.8}}
\put (10,0){\circle*{1}}
\end{picture} }
&$\emptyset$ &{$\tableau[b,s]{\\  }$} &{$\tableau[b,s]{1}$}
&${\begin{array}{c} |Y_{1,1}| \\ \quad \\ \quad \\ \end{array}}$
&${\begin{array}{c} y_{11} \\ \quad \\ \quad \\ \end{array}}$
&${\begin{array}{c} e(T) \\ \quad \\ \quad  \\ \end{array}}$ \\
\hline\hline

\end{tabular}

\begin{tabular}{|c|c|c|c|c|c|c|c|} \hline\hline
\multicolumn{8}{|c|}{Generators of $SL_4$ Tensor Product Algebra}
\\ \hline\hline
{No.} &{$\ $ BZ Diagrams $\ $} &\ $D^t\ $  &\ $E^t\ $  &\ $F^t\ $
&$\Delta_{(D,E,F), M(T)}$ &$e(T)$ &${\mathcal E}(T)$ \\
\hline\hline &&&&&&& \\

${\begin{array}{c} 4 \\ \quad \\ \quad \\ \end{array}}$
& {\begin{picture}(10,7.5)
\thicklines
\put (0,0){\line(1,0){10}}
\put (0,0){\line(3,5){5}}
\put (4,0){\line(3,5){3}}
\put (8,0){\line(3,5){1}}
\put (2,0){\line(-3,5){1}}
\put (6,0){\line(-3,5){3}}
\put (10,0){\line(-3,5){5}}
\put (2.1,3.464){\line(1,0){5.8}}
\put (4.1,6.928){\line(1,0){1.8}}
\thinlines
\put (3,5.196){\line(1,0){4}}
\put (3,5.196){\circle*{1}}
\put (7,5.196){\circle*{1}}
\end{picture} }
&{$\tableau[b,s]{\\  \\ }$} &$\emptyset$ &{$\tableau[b,s]{\\ \\
}$} &${\begin{array}{c} |X_{2,2}|  \\ \quad \\ \quad \\
\end{array}}$
&${\begin{array}{c} 1 \\ \quad \\ \quad
\\ \end{array}}$
&${\begin{array}{c} x_{11}x_{22}e(T) \\ \quad \\ \quad
\\ \end{array}}$\\
\hline\hline    &&&&&&& \\

${\begin{array}{c} 5 \\ \quad \\ \quad \\ \end{array}}$
& {\begin{picture}(10,7.5)
\thicklines
\put (0,0){\line(1,0){10}}
\put (0,0){\line(3,5){5}}
\put (4,0){\line(3,5){3}}
\put (8,0){\line(3,5){1}}
\put (2,0){\line(-3,5){1}}
\put (6,0){\line(-3,5){3}}
\put (10,0){\line(-3,5){5}}
\put (2.1,3.464){\line(1,0){5.8}}
\put (4.1,6.928){\line(1,0){1.8}}
\thinlines
\put (1,1.732){\line(1,0){8}}
\put (1,1.732){\circle*{1}}
\put (5,1.732){\circle*{1}}
\put (9,1.732){\circle*{1}}
\end{picture} }
&{$\tableau[b,s]{\\  }$} &$\emptyset$ &{$\tableau[b,s]{\\ }$}
&${\begin{array}{c} |X_{1,1}| \\ \quad \\ \quad \\ \end{array}}$
&${\begin{array}{c} 1 \\ \quad \\ \quad \\
\end{array}}$
&${\begin{array}{c} x_{11}e(T) \\ \quad \\ \quad \\
\end{array}}$\\
\hline\hline &&&&&&& \\

${\begin{array}{c} 6 \\ \quad \\ \quad \\ \end{array}}$ &
{\begin{picture}(10,7.5) \thicklines \put (0,0){\line(1,0){10}}
\put (0,0){\line(3,5){5}} \put (4,0){\line(3,5){3}} \put
(8,0){\line(3,5){1}} \put (2,0){\line(-3,5){1}} \put
(6,0){\line(-3,5){3}} \put (10,0){\line(-3,5){5}} \put
(2.1,3.464){\line(1,0){5.8}} \put (4.1,6.928){\line(1,0){1.8}}
\thinlines \put (4,0){\line(-3,5){2}} \put (4,0){\circle*{1}} \put
(2,3.464){\circle*{1}}
\end{picture} }
&{$\tableau[b,s]{ \\ \\ }$} &{$\tableau[b,s]{\\ \\ }$}
&{$\tableau[b,s]{ \\ \\ 1 \\ 2 }$} &${\begin{array}{c}
|X_{4,2}Y_{4,2}| \\ \quad \\ \quad \\ \end{array}}$
&${\begin{array}{c} y_{31}y_{42} \\ \quad \\ \quad \\
\end{array}}$
&${\begin{array}{c} x_{11}x_{22} e(T) \\ \quad \\ \quad \\
\end{array}}$
\\ \hline\hline &&&&&&& \\

${\begin{array}{c} 7 \\ \quad \\ \quad \\ \end{array}}$
& {\begin{picture}(10,7.5)
\thicklines
\put (0,0){\line(1,0){10}}
\put (0,0){\line(3,5){5}}
\put (4,0){\line(3,5){3}}
\put (8,0){\line(3,5){1}}
\put (2,0){\line(-3,5){1}}
\put (6,0){\line(-3,5){3}}
\put (10,0){\line(-3,5){5}}
\put (2.1,3.464){\line(1,0){5.8}}
\put (4.1,6.928){\line(1,0){1.8}}
\thinlines
\put (8,0){\line(-3,5){4}}
\put (8,0){\circle*{1}}
\put (6,3.464){\circle*{1}}
\put (4,6.928){\circle*{1}}
\end{picture} }
&{$\tableau[b,s]{\\ \\ \\ }$} &{$\tableau[b,s]{\\  }$} &{$\tableau[b,s]{ \\ \\ \\ 1 }$} &${\begin{array}{c} |X_{4,3}Y_{4,1}| \\ \quad \\ \quad \\ \end{array}}$
&${\begin{array}{c} y_{41} \\ \quad \\ \quad \\ \end{array}}$
&${\begin{array}{c} x_{11}x_{22}x_{33} e(T) \\ \quad \\ \quad \\
\end{array}}$
\\ \hline\hline &&&&&&& \\

${\begin{array}{c} 8 \\ \quad \\ \quad \\ \end{array}}$
& {\begin{picture}(10,7.5)
\thicklines
\put (0,0){\line(1,0){10}}
\put (0,0){\line(3,5){5}}
\put (4,0){\line(3,5){3}}
\put (8,0){\line(3,5){1}}
\put (2,0){\line(-3,5){1}}
\put (6,0){\line(-3,5){3}}
\put (10,0){\line(-3,5){5}}
\put (2.1,3.464){\line(1,0){5.8}}
\put (4.1,6.928){\line(1,0){1.8}}
\thinlines
\put (6,0){\line(3,5){2}}
\put (6,0){\circle*{1}}
\put (8,3.464){\circle*{1}}
\end{picture} }
&$\emptyset$ &{$\tableau[b,s]{\\ \\ }$} &{$\tableau[b,s]{1 \\ 2
}$} &${\begin{array}{c} |Y_{2,2}| \\ \quad \\ \quad \\
\end{array}}$
&${\begin{array}{c} y_{11}y_{22} \\ \quad \\ \quad \\
\end{array}}$
&${\begin{array}{c} e(T) \\ \quad \\ \quad \\
\end{array}}$
 \\ \hline\hline
\end{tabular}

\vskip 20pt

\begin{tabular}{|c|c|c|c|c|c|c|c|} \hline\hline
\multicolumn{8}{|c|}{Generators of $SL_4$ Tensor Product Algebra}
\\ \hline\hline
{No.} &{$\ $ BZ Diagrams $\ $} &\ $D^t\ $  &\ $E^t\ $  &\ $F^t\ $
&$\Delta_{(D,E,F), M(T)}$ &$e(T)$ &${\mathcal E}(T)$ \\
\hline\hline &&&&&&& \\

${\begin{array}{c} 9 \\ \quad \\ \quad \\ \end{array}}$
& {\begin{picture}(10,7.5)
\thicklines
\put (0,0){\line(1,0){10}}
\put (0,0){\line(3,5){5}}
\put (4,0){\line(3,5){3}}
\put (8,0){\line(3,5){1}}
\put (2,0){\line(-3,5){1}}
\put (6,0){\line(-3,5){3}}
\put (10,0){\line(-3,5){5}}
\put (2.1,3.464){\line(1,0){5.8}}
\put (4.1,6.928){\line(1,0){1.8}}
\thinlines
\put (2,0){\line(3,5){4}}
\put (4,3.464){\circle*{1}}
\put (2,0){\circle*{1}}
\put (6,6.928){\circle*{1}}
\end{picture} }
&$\emptyset$ &{$\tableau[b,s]{\\ \\ \\ }$} &{$\tableau[b,s]{1 \\ 2 \\
3
}$} &${\begin{array}{c} |Y_{3,3}| \\ \quad \\ \quad \\
\end{array}}$
&${\begin{array}{c} y_{11}y_{22}y_{33} \\ \quad \\ \quad \\
\end{array}}$
&${\begin{array}{c} e(T) \\ \quad \\ \quad \\
\end{array}}$
 \\ \hline\hline   &&&&&&& \\

${\begin{array}{c} 10 \\ \quad \\ \quad \\ \end{array}}$
& {\begin{picture}(10,7.5)
\thicklines
\put (0,0){\line(1,0){10}}
\put (0,0){\line(3,5){5}}
\put (4,0){\line(3,5){3}}
\put (8,0){\line(3,5){1}}
\put (2,0){\line(-3,5){1}}
\put (6,0){\line(-3,5){3}}
\put (10,0){\line(-3,5){5}}
\put (2.1,3.464){\line(1,0){5.8}}
\put (4.1,6.928){\line(1,0){1.8}}
\thinlines
\put (2,0){\line(3,5){3}}
\put (4,6.928){\line(3,-5){1}}
\put (7,5.196){\line(-1,0){2}}
\put (2,0){\circle*{1}}
\put (4,3.464){\circle*{1}}
\put (4,6.928){\circle*{1}}
\put (7,5.196){\circle*{1}}
\end{picture} }
&{$\tableau[b,s]{\\ \\ \\ }$} &{$\tableau[b,s]{\\ \\ \\ }$} &{$\tableau[b,s]{ &1 \\ &2 \\ \\ 3 }$}
&${\begin{array}{c} {\left| \begin{array}{cc} X_{4,3} & Y_{4,3} \\ 0 &Y_{23}\end{array} \right|} \\ \quad \\ \quad \\ \end{array}}$
&${\begin{array}{c} y_{11}y_{22}y_{43} \\ \quad \\ \quad \\ \end{array}}$
&${\begin{array}{c} x_{11}x_{22}x_{33} e(T) \\ \quad \\ \quad \\
\end{array}}$
 \\ \hline\hline &&&&&&& \\

${\begin{array}{c} 11 \\ \quad \\ \quad \\ \end{array}}$ &
{\begin{picture}(10,7.5) \thicklines \put (0,0){\line(1,0){10}}
\put (0,0){\line(3,5){5}} \put (4,0){\line(3,5){3}} \put
(8,0){\line(3,5){1}} \put (2,0){\line(-3,5){1}} \put
(6,0){\line(-3,5){3}} \put (10,0){\line(-3,5){5}} \put
(2.1,3.464){\line(1,0){5.8}} \put (4.1,6.928){\line(1,0){1.8}}
\thinlines \put (8,0){\line(-3,5){3}} \put
(3,5.196){\line(1,0){2}} \put (6,6.928){\line(-3,-5){1}} \put
(8,0){\circle*{1}} \put (6,3.464){\circle*{1}} \put
(3,5.196){\circle*{1}} \put (6,6.928){\circle*{1}}
\end{picture} }
&{$\tableau[b,s]{\\ \\ }$} &{$\tableau[b,s]{\\  }$}
&{$\tableau[b,s]{ \\ \\ 1}$} &${\begin{array}{c} |X_{3,2}Y_{3,1}|
\\ \quad \\ \quad \\ \end{array}}$ &${\begin{array}{c} y_{31} \\
\quad \\ \quad \\ \end{array}}$
&${\begin{array}{c} x_{11}x_{22} e(T) \\ \quad \\ \quad \\
\end{array}}$
 \\ \hline\hline   &&&&&&& \\

${\begin{array}{c} 12 \\ \quad \\ \quad \\ \end{array}}$ &
{\begin{picture}(10,7.5) \thicklines \put (0,0){\line(1,0){10}}
\put (0,0){\line(3,5){5}} \put (4,0){\line(3,5){3}} \put
(8,0){\line(3,5){1}} \put (2,0){\line(-3,5){1}} \put
(6,0){\line(-3,5){3}} \put (10,0){\line(-3,5){5}} \put
(2.1,3.464){\line(1,0){5.8}} \put (4.1,6.928){\line(1,0){1.8}}
\thinlines \put (2,0){\line(3,5){1}} \put
(2,3.464){\line(3,-5){1}} \put (9,1.732){\line(-1,0){6}} \put
(2,0){\circle*{1}} \put (2,3.464){\circle*{1}} \put
(5,1.732){\circle*{1}} \put (9,1.732){\circle*{1}}
\end{picture} }
&{$\tableau[b,s]{\\ \\ }$} &{$\tableau[b,s]{\\ \\ \\ }$}
&{$\tableau[b,s]{ &1 \\ \\ 2\\ 3 }$} &${\begin{array}{c} {\left|
\begin{array}{cc} X_{4,2} & Y_{4,3} \\ 0 &Y_{13}\end{array}
\right|} \\ \quad \\ \quad \\ \end{array}}$
 &${\begin{array}{c} y_{11}y_{32}y_{43} \\ \quad \\ \quad \\ \end{array}}$
&${\begin{array}{c} x_{11}x_{22} e(T) \\ \quad \\ \quad \\
\end{array}}$
 \\ \hline\hline &&&&&&& \\

${\begin{array}{c} 13 \\ \quad \\ \quad \\ \end{array}}$ &
{\begin{picture}(10,7.5) \thicklines \put (0,0){\line(1,0){10}}
\put (0,0){\line(3,5){5}} \put (4,0){\line(3,5){3}} \put
(8,0){\line(3,5){1}} \put (2,0){\line(-3,5){1}} \put
(6,0){\line(-3,5){3}} \put (10,0){\line(-3,5){5}} \put
(2.1,3.464){\line(1,0){5.8}} \put (4.1,6.928){\line(1,0){1.8}}
\thinlines \put (1,1.732){\line(1,0){2}} \put
(4,0){\line(-3,5){1}} \put (6,6.928){\line(-3,-5){3}} \put
(4,3.464){\circle*{1}} \put (1,1.732){\circle*{1}} \put
(4,0){\circle*{1}} \put (6,6.928){\circle*{1}}
\end{picture} }
&{$\tableau[b,s]{\\  }$} &{$\tableau[b,s]{\\  \\ }$} &{$\tableau[b,s]{ \\ 1\\ 2}$} &${\begin{array}{c} |X_{3,1}Y_{3,2}| \\ \quad \\ \quad \\ \end{array}}$
&${\begin{array}{c} y_{21}y_{32} \\ \quad \\ \quad \\ \end{array}}$
&${\begin{array}{c} x_{11}e(T) \\ \quad \\ \quad \\
\end{array}}$
 \\ \hline\hline

\end{tabular}

\begin{tabular}{|c|c|c|c|c|c|c|c|} \hline\hline
\multicolumn{8}{|c|}{Generators of $SL_4$ Tensor Product Algebra}
\\ \hline\hline
{No.} &{BZ Diagrams} &$D^t$  &$E^t$  &$F^t$
&$\Delta_{(D,E,F), M(T)}$ &$e(T)$ &${\mathcal E}(T)$ \\
\hline\hline &&&&&&& \\

${\begin{array}{c} 14 \\ \quad \\ \quad \\ \end{array}}$
& {\begin{picture}(10,7.5)
\thicklines
\put (0,0){\line(1,0){10}}
\put (0,0){\line(3,5){5}}
\put (4,0){\line(3,5){3}}
\put (8,0){\line(3,5){1}}
\put (2,0){\line(-3,5){1}}
\put (6,0){\line(-3,5){3}}
\put (10,0){\line(-3,5){5}}
\put (2.1,3.464){\line(1,0){5.8}}
\put (4.1,6.928){\line(1,0){1.8}}
\thinlines
\put (6,0){\line(3,5){1}}
\put (4,6.928){\line(3,-5){3}}
\put (9,1.732){\line(-1,0){2}}
\put (6,3.464){\circle*{1}}
\put (6,0){\circle*{1}}
\put (4,6.928){\circle*{1}}
\put (9,1.732){\circle*{1}}
\end{picture} }
&{$\tableau[b,s]{\\ \\ \\ }$} &{$\tableau[b,s]{\\ \\ }$} &{$\tableau[b,s]{ &1 \\ \\ \\ 2}$}
&${\begin{array}{c}  {\left| \begin{array}{cc} X_{4,3} & Y_{4,2} \\ 0 &Y_{1,2}\end{array} \right|} \\ \quad \\ \quad \\ \end{array}}$
&${\begin{array}{c} y_{11}y_{42} \\ \quad \\ \quad \\ \end{array}}$
&${\begin{array}{c} x_{11}x_{22}x_{33} e(T) \\ \quad \\ \quad \\
\end{array}}$
 \\ \hline\hline   &&&&&&& \\

${\begin{array}{c} 15 \\ \quad \\ \quad \\ \end{array}}$
& {\begin{picture}(10,7.5)
\thicklines
\put (0,0){\line(1,0){10}}
\put (0,0){\line(3,5){5}}
\put (4,0){\line(3,5){3}}
\put (8,0){\line(3,5){1}}
\put (2,0){\line(-3,5){1}}
\put (6,0){\line(-3,5){3}}
\put (10,0){\line(-3,5){5}}
\put (2.1,3.464){\line(1,0){5.8}}
\put (4.1,6.928){\line(1,0){1.8}}
\thinlines
\put (8,0){\line(-3,5){1}}
\put (8,3.464){\line(-3,-5){1}}
\put (1,1.732){\line(1,0){6}}
\put (8,0){\circle*{1}}
\put (8,3.464){\circle*{1}}
\put (5,1.732){\circle*{1}}
\put (1,1.732){\circle*{1}}
\end{picture} }
&{$\tableau[b,s]{\\  }$} &{$\tableau[b,s]{\\  }$} &{$\tableau[b,s]{ \\ 1 }$} &${\begin{array}{c} |X_{2,1}Y_{2,1}| \\ \quad \\ \quad \\ \end{array}}$
&${\begin{array}{c} y_{21} \\ \quad \\ \quad \\ \end{array}}$
&${\begin{array}{c} x_{11}e(T) \\ \quad \\ \quad \\
\end{array}}$
 \\ \hline\hline &&&&&&& \\

${\begin{array}{c} 16 \\ \quad \\ \quad \\ \end{array}}$
& {\begin{picture}(10,7.5)
\thicklines
\put (0,0){\line(1,0){10}}
\put (0,0){\line(3,5){5}}
\put (4,0){\line(3,5){3}}
\put (8,0){\line(3,5){1}}
\put (2,0){\line(-3,5){1}}
\put (6,0){\line(-3,5){3}}
\put (10,0){\line(-3,5){5}}
\put (2.1,3.464){\line(1,0){5.8}}
\put (4.1,6.928){\line(1,0){1.8}}
\thinlines
\put (6,0){\line(3,5){1}}
\put (6,6.928){\line(-3,-5){1}}
\put (3,5.196){\line(1,0){2}}
\put (9,1.732){\line(-1,0){2}}
\put (5,5.196){\line(3,-5){2}}
\put (6,0){\circle*{1}}
\put (6,3.464){\circle*{1}}
\put (6,6.928){\circle*{1}}
\put (3,5.196){\circle*{1}}
\put (9,1.732){\circle*{1}}
\end{picture} }
&{$\tableau[b,s]{\\ \\ }$} &{$\tableau[b,s]{\\ \\ }$} &{$\tableau[b,s]{ &1 \\ \\ 2 }$}
&${\begin{array}{c} {\left| \begin{array}{cc} X_{3,2} & Y_{3,2} \\ 0 &Y_{1,2}\end{array} \right|}  \\ \quad \\ \quad \\ \end{array}}$
&${\begin{array}{c} y_{11}y_{32} \\ \quad \\ \quad \\ \end{array}}$
&${\begin{array}{c} x_{11}x_{22} e(T) \\ \quad \\ \quad \\
\end{array}}$
 \\ \hline\hline   &&&&&&& \\

${\begin{array}{c} 17 \\ \quad \\ \quad \\ \end{array}}$
& {\begin{picture}(10,7.5)
\thicklines
\put (0,0){\line(1,0){10}}
\put (0,0){\line(3,5){5}}
\put (4,0){\line(3,5){3}}
\put (8,0){\line(3,5){1}}
\put (2,0){\line(-3,5){1}}
\put (6,0){\line(-3,5){3}}
\put (10,0){\line(-3,5){5}}
\put (2.1,3.464){\line(1,0){5.8}}
\put (4.1,6.928){\line(1,0){1.8}}
\thinlines
\put (1,1.732){\line(1,0){2}}
\put (4,0){\line(-3,5){1}}
\put (3,1.732){\line(3,5){2}}
\put (7,5.196){\line(-1,0){2}}
\put (4,6.928){\line(3,-5){1}}
\put (1,1.732){\circle*{1}}
\put (4,3.464){\circle*{1}}
\put (4,6.928){\circle*{1}}
\put (4,0){\circle*{1}}
\put (7,5.196){\circle*{1}}
\end{picture} }
&{$\tableau[b,s]{ & \\ \\ \\ }$} &{$\tableau[b,s]{\\ \\ }$} &{$\tableau[b,s]{ & \\ &1 \\ \\ 2 }$}
&${\begin{array}{c}  {\left| \begin{array}{ccc} X_{4,3} &0 &Y_{4,2} \\ 0 &X_{2,1} &Y_{2,2}\end{array} \right|}   \\ \quad \\ \quad \\ \end{array}}$
&${\begin{array}{c} y_{21}y_{42} \\ \quad \\ \quad \\ \end{array}}$
&${\begin{array}{c} x_{11}^2x_{22}x_{33} e(T) \\ \quad \\ \quad \\
\end{array}}$
 \\ \hline\hline &&&&&&& \\

${\begin{array}{c} 18 \\ \quad \\ \quad \\ \end{array}}$
& {\begin{picture}(10,7.5)
\thicklines
\put (0,0){\line(1,0){10}}
\put (0,0){\line(3,5){5}}
\put (4,0){\line(3,5){3}}
\put (8,0){\line(3,5){1}}
\put (2,0){\line(-3,5){1}}
\put (6,0){\line(-3,5){3}}
\put (10,0){\line(-3,5){5}}
\put (2.1,3.464){\line(1,0){5.8}}
\put (4.1,6.928){\line(1,0){1.8}}
\thinlines
\put (2,0){\line(3,5){1}}
\put (8,0){\line(-3,5){1}}
\put (3,1.732){\line(1,0){4}}
\put (2,3.464){\line(3,-5){1}}
\put (8,3.464){\line(-3,-5){1}}
\put (5,1.732){\circle*{1}}
\put (2,0){\circle*{1}}
\put (8,0){\circle*{1}}
\put (2,3.464){\circle*{1}}
\put (8,3.464){\circle*{1}}
\end{picture} }
&{$\tableau[b,s]{\\ \\  }$} &{$\tableau[b,s]{ & \\ \\ \\ }$} &{$\tableau[b,s]{ &1 \\ &2 \\ 1 \\ 3 }$}
&${\begin{array}{c} {\left| \begin{array}{ccc} X_{4,2} &Y_{4,3} &Y_{4,1} \\ 0 &Y_{2,3} &0 \end{array} \right|}  \\ \quad \\ \quad \\ \end{array}}$
&${\begin{array}{c} y_{11}y_{22}y_{31}y_{43} \\ \quad \\ \quad \\ \end{array}}$
&${\begin{array}{c} x_{11}x_{22} e(T) \\ \quad \\ \quad \\
\end{array}}$
 \\ \hline\hline
\end{tabular}

\end{document}